\documentclass[review,3p,12]{elsarticle}

\usepackage{amssymb}
\usepackage{amsmath}
\usepackage{float}
\usepackage[numbers]{natbib}
\usepackage{comment}
\usepackage{adjustbox}
\usepackage{afterpage}
\usepackage{graphicx}
\usepackage{subcaption}
\usepackage{booktabs}
\usepackage{soul}
\usepackage{xcolor}
\usepackage{caption}
\usepackage{tabularx}
\usepackage{array}
\usepackage{layout}
\usepackage[colorlinks=true, linkcolor=blue, citecolor=blue, filecolor=blue, urlcolor=blue]{hyperref}

\biboptions{sort&compress}
\renewcommand{\vec}[1]{\ensuremath\boldsymbol{#1}}

\usepackage{lineno}

\journal{Applied Numerical Mathematics}

\begin{document}

\begin{frontmatter}

\title{A Modified Hermite Radial Basis Function for Accurate Interpolation}

\author[1]{Amirhossein Fashamiha}
\author[1]{David Salac\corref{cor1}}
\ead{davidsal@buffalo.edu}
\cortext[cor1]{Corresponding author}

\affiliation[1]{organization={Department of Mechanical and Aerospace Engineering, University at Buffalo},
            city={Buffalo},
            state={NY},
            postcode={14260-4400},
            country={United States}}

\begin{abstract}
Accurate interpolation of functions and derivatives is crucial in solving partial differential equations (PDEs). The Radial Basis Function (RBF) method has become an extremely popular and robust approach for interpolation on scattered data. Hermite Radial Basis Function (HRBF) methods are an extension of the RBF and improve the overall accuracy by incorporating both function and derivative information. Infinitely smooth kernels, such as the Gaussian, use a shape-parameter to describe the width of support and are widely used due to their excellent approximation accuracy and ability to capture fine-scale details. Unfortunately, the use of infinitely smooth kernels suffers from ill-conditioning at low to moderate shape parameters, which affects the accuracy. This work proposes a Modified HRBF (MHRBF) method that introduces an additional polynomial term to balance kernel behavior, improving accuracy while maintaining or lowering computational cost. Using standard double-precision mathematics, the results indicate that compared to the HRBF method, the MHRBF method achieves lower error for all values of the shape parameter and domain size. The MHRBF is also able to achieve low errors at a lower computational cost as compared to the standard HRBF method.
\end{abstract}

\begin{keyword}
Hermite Radial Basis Functions \sep Numerical Stability \sep High-Order Interpolation \sep Computational Cost


\end{keyword}

\end{frontmatter}


\section{Introduction}

Radial Basis Functions (RBFs) are a powerful approach for multivariate interpolation, providing a robust framework for accurately reconstructing smooth functions from scattered data. Although the term `Radial Basis Function' was formally introduced by Dyn and Levin in 1983~\cite{Dyn1983}, key functions now classified as RBFs had already been studied in the 1970s. Hardy (1971) pioneered the multiquadric (MQ) method, using sums of radial quadric functions to interpolate scattered geospatial data and reconstruct complex topography and irregular surfaces~\cite{Hardy1971}. Later, Duchon (1977) established the theoretical foundation of RBFs through the development of the thin plate splines (TPS) method~\cite{duchon1977splines}. Franke's well-known benchmark studies~\cite{franke1979critical, franke1982scattered} further highlighted the practical strengths of RBF methods, showing that the MQ RBF often outperformed other interpolation approaches, particularly for scattered or irregular data.

Achieving high-order approximations of functions and their derivatives is critical for solving partial differential equations (PDEs)~\cite{kansa1990multiquadrics, mai2003approximation, kansa2004volumetric, li2013localized, dehghan2014numerical}. Classical mesh-based methods such as finite difference, finite element, and finite volume rely on structured grids, mesh connectivity, or surface integrations, which complicate their use on complex or evolving geometries. In contrast, RBF methods are inherently mesh-free, enabling high-order accurate approximations directly on scattered nodes. This flexibility makes RBFs particularly well-suited for solving problems on irregular domains. The application of RBFs to PDEs attracted major notice following Kansa's work~\cite{kansa1990multiquadrics}, which introduced the use of global RBFs in collocation methods. However, Kansa's method is known to suffer from stability and ill-conditioning issues as the number of nodes increases. Later, Fasshauer addressed the ill-conditioning of Kansa's RBF collocation method by introducing sparse system matrices and symmetric formulations~\cite{fasshauer1996solving}.

Hermite interpolation, introduced by Charles Hermite in 1878, evaluates a function and its derivatives at specific points using predefined nodal values~\cite{davis1975interpolation}. Building on this, Birkhoff proposed Hermite-Birkhoff interpolation~\cite{birkhoff1906general}, enabling the interpolation of higher-order derivatives based solely on prescribed nodal values of those derivatives. Wu later extended RBFs to the Hermite-Birkhoff framework, creating the Hermite RBF (HRBF) approach, which combines RBFs with their derivatives~\cite{zongmin1992hermite}. This HRBF method enables both function and derivative interpolation, making it well-suited to handle gradient and higher-order derivative constraints in applications.

In this paper, we propose a Modified Hermite Radial Basis Function (MHRBF) method that introduces spatially varying coefficients into the kernel, enabling improved interpolation accuracy without increasing computational cost. The goal is to increase the stability and accuracy of the method using double-precision mathematics without resorting to methods such as RBF-QR or extended precision mathematics, both of which are computationally expensive. Unlike conventional HRBF methods, which rely solely on the RBF and its derivatives, the proposed approach modifies the kernel to explicitly control localized behavior through multiplicative coefficients. By modifying the influence of the kernel near interpolation points and balancing it across the domain, this modification enhances the ability of the method to capture fine details while maintaining or even reducing computational expense. Furthermore, the new formulation reduces the dependency on higher-order derivatives by limiting their contribution to one order lower than in conventional HRBF methods. In general, the MHRBF method bridges the gap between accuracy and computational cost, offering a robust approach to interpolation problems.

\section{From Radial Basis Functions to Modified Hermite Radial Basis Functions}

RBFs are a class of functions used for interpolating multivariate data in a mesh-free manner. For a set of scattered points, an RBF interpolant can be expressed by a linear combination of an arbitrary radial function in a $d$-dimensional space:
\begin{equation}\label{rbf}
    s(\vec{x})=\sum_{i=1}^N w_i \phi(\|\vec{x}-\vec{x}_i\|).
\end{equation}
Here, $w_i$ is a set of weights, $\phi(\|\vec{x}-\vec{x}_i\|)$ is the RBF kernel applied over all the given scattered points $N$, and the norm $\| . \|$ is the 2-norm. In standard RBFs, the only constraint for function representation is $s(\vec{x}_p)=f_p$, where $(\vec{x}_p,f_p)$ for $p\in[1,N]$ represents the set of given data. The RBF weights, $w_i$, are calculated by enforcing this constraint, which leads to a linear system, $\vec{A}\vec{w}=\vec{f}$, that must be solved for a given data set:
\begin{equation}\label{rbfmatrix}
\underbrace{
\begin{bmatrix}
\phi(\|\vec{x}_{1}-\vec{x}_{1}\|) & \phi(\|\vec{x}_{1}-\vec{x}_{2}\|) & \dots & \phi(\|\vec{x}_{1}-\vec{x}_{N}\|)\\
\phi(\|\vec{x}_{2}-\vec{x}_{1}\|) & \phi(\|\vec{x}_{2}-\vec{x}_{2}\|) & \dots & \phi(\|\vec{x}_{2}-\vec{x}_{N}\|)\\
\vdots & \vdots & \vdots & \vdots\\
\phi(\|\vec{x}_{N}-\vec{x}_{1}\|) & \phi(\|\vec{x}_{N}-\vec{x}_{2}\|) & \dots & \phi(\|\vec{x}_{N}-\vec{x}_{N}\|)\\
\end{bmatrix}}_{\text{$\vec{A}$}}
\underbrace{
\begin{bmatrix}
w_{1}\\
w_{2}\\
\vdots\\
w_{N}\\
\end{bmatrix}}_{\text{$\vec{w}$}}
=
\underbrace{
\begin{bmatrix}
f(\vec{x}_{1})\\
f(\vec{x}_{2})\\
\vdots\\
f(\vec{x}_{N})\\
\end{bmatrix}}_{\text{$\vec{f}$}}.
\end{equation}
Sample RBF kernels are shown in~\autoref{kernels}, where $r=\|\vec{x}-\vec{x}_i\|$ represents the distance between two points, and $\varepsilon$ is the shape parameter.
\begin{table}[h]
\centering
\caption{Classes of Basis Functions and Their Radial Functions}
\resizebox{0.9\textwidth}{!}{
\begin{tabular}{>{\raggedright\arraybackslash}p{0.45\textwidth} >{\raggedright\arraybackslash}p{0.5\textwidth}}
\toprule
Class of Basis Function & Radial Function $\phi(r)$ \\
\midrule
Gaussian (GA)           & $e^{-(\varepsilon r)^2}$ \\
Polyharmonic Spline (PHS) & $r^{2k-1}$ or $r^{2k}\log r$, $k \in \mathbb{N}$ \\
Multiquadric (MQ)       & $\sqrt{1+(\varepsilon r)^2}$\\
Inverse Multiquadric (IMQ) & $\frac{1}{\sqrt{1+(\varepsilon r)^2}}$ \\
Inverse Quadric (IQ)    & $\frac{1}{1+(\varepsilon r)^2}$ \\
Bessel (BE)             & $\frac{J_{p/2-1}(\varepsilon r)}{(\varepsilon r)^{p/2-1}}$ \\
\bottomrule
\end{tabular}
}
\label{kernels}
\end{table}

In \autoref{kernels}, the shape parameter ($\varepsilon$) for infinitely smooth RBFs is usually related to a representative grid spacing, $h$, and influences both the flatness of the function and the accuracy of the interpolant. Infinitely smooth RBFs are widely used due to their superior approximation accuracy and ability to capture fine-scale features~\cite{larsson2003numerical}. However, their sensitivity to $\varepsilon$ introduces major challenges. As $\varepsilon$ decreases, the kernel flattens, leading to ill-conditioning and stagnation errors, where accuracy ceases to improve despite increasing node density. In contrast, larger $\varepsilon$ values localize the kernel, improving conditioning but increasing interpolation error due to underfitting. One common alternative is piecewise smooth RBFs such as PHS, which do not require a shape parameter and are less prone to ill-conditioning. However, these functions lack the exponential convergence rates of infinitely smooth RBFs, making them unsuitable for high accuracy applications that benefit from smooth kernels like GA and MQ~\cite{fornberg2007runge}.

Several advanced methods have been developed to address the stability issues in RBF interpolation, including RBF-CP (Contour Padé)~\cite{fornberg2004stable}, RBF-QR (Rational Quadrature)~\cite{driscoll2002interpolation, larsson2013stable}, RBF-GA (Gaussian Approximation)~\cite{fornberg2013stable}, or through the use of extended precision~\cite{flyer2016role}. However, the implementation of these methods can be computationally intensive. In addition to these techniques, RBF-FD (Finite Differences) \cite{bayona20153, chandhini2007local, chinchapatnam2009compact, roque2011local} offers a method that combines RBF with finite difference approaches to solve PDEs in irregular domains, providing greater flexibility and reliability in complex geometries.

\subsection{RBF Interpolation with Polynomial Augmentation}

The RBF interpolant in \autoref{rbf} is typically augmented with a polynomial basis to ensure that the interpolation matrix (\autoref{rbfmatrix}) remains non-singular, particularly when working with conditionally positive definite RBFs, such as PHS. Polynomial augmentation not only stabilizes the interpolation matrix but also improves the accuracy of derivative approximations by reducing oscillations, which is especially beneficial for achieving smoother solutions. For a more comprehensive discussion of the role of polynomials in RBF interpolation, refer to \cite{flyer2016role}. With polynomial augmentation, the RBF interpolant in \autoref{rbf} can be re-written as a combination of a radial basis function term and a polynomial term:
\begin{equation}\label{augment}
s(\vec{x})=\sum^{N}_{i=1}w_{i}\phi(\|\vec{x}-\vec{x}_{i}\|)+\sum^{M}_{k=1}\lambda_{k}p_k(\vec{x}),
\end{equation}
where the set $p_k(\vec{x})$ are the basis function for polynomials of degree up to $l$ in $d$-dimensional space, of which there are $M=\binom{l+d}{d}$, and $\lambda_{k}$ represents the polynomial weights. The addition of a constraint
\begin{equation} \label{eq:polyAugConstraint}
    \sum^{N}_{i=1}w_i p_k(\vec{x}_i)=0 \textnormal{ for } k\in[1,M]
\end{equation}
results in the following linear system:
\begin{equation}\label{RBF_augmentedmatrix}
\begin{bmatrix}
\vec{A} & \vec{P}\\
\vec{P}^T & \vec{0}\\
\end{bmatrix}
\begin{bmatrix}
\vec{w}\\
\vec{\lambda}\\
\end{bmatrix}
=
\begin{bmatrix}
\vec{f}\\
\vec{0}\\
\end{bmatrix}
\end{equation}
where $\vec{P}$ is the matrix of polynomial terms given by $(\vec{P})_{ij}=p_{j}(\vec{x}_i)$, with $i\in[1,N]$ and $j\in[1,M]$.

\subsection{Hermite Radial Basis Functions}
The HRBF interpolant extends the standard RBF by incorporating constraints on both the function values and gradients at each interpolation point $\vec{x}_p$, ensuring gradient continuity; an essential feature for derivative-based problems like PDEs~\cite{lehto2017radial}. Specifically, the HRBF satisfies $s(\vec{x}_p)=f_p$ and $\nabla s(\vec{x}_p)=\vec{g}_p$, where $f_p$ and $\vec{g}_p$ are the known function values and gradients, respectively, at $\vec{x}_p$.  The formulation for HRBF can be expressed as:
\begin{equation}\label{hermite-first}
s(\vec{x})=\sum^{N}_{i=1}\left(w_{i}\phi(\|\vec{x}- \vec{x}_{i}\|)+\vec{b}_{i}\cdot\nabla\phi(\|\vec{x}-\vec{x}_{i}\|)\right)+\sum^{M}_{k=1}\lambda_{k}p_k(\vec{x}).
\end{equation}
The second term incorporates gradient information, with $\nabla$ as the gradient operator acting on $\phi$, and
$\vec{b}_{i} = \begin{bmatrix} \alpha_{i} & \beta_{i} & \gamma_{i} \end{bmatrix}^T$ (in 3D) representing the weight vector associated with the gradients at the interpolation points $\vec{x}_{i}$.

In two dimensions, computing the HRBF weights requires solving the following linear system:
\begin{equation}\label{augmentedmatrix}
\begin{bmatrix}
    \vec{A} & \vec{A}_x & \vec{A}_y & \vec{P}\\
    \vec{A}_x & \vec{A}_{xx} & \vec{A}_{yy} & \vec{P}_x\\
    \vec{A}_y & \vec{A}_{xy} & \vec{A}_{yy} & \vec{P}_y\\
    \vec{P}^T & \vec{P}^T_x & \vec{P}^T_y & \vec{0}
\end{bmatrix}
\begin{bmatrix}
   \vec{w} \\
   \vec{\alpha}\\
   \vec{\beta}\\
   \vec{\lambda}
\end{bmatrix}
=
\begin{bmatrix}
    \vec{f}\\
    \vec{f}_x\\
    \vec{f}_y\\
    \vec{0}
\end{bmatrix},
\end{equation}
where $\vec{A}$ is the matrix of the RBF kernel and $\vec{P}$ is the matrix of polynomial terms, with subscripts representing the derivatives of the kernel or polynomial, as appropriate. Micchelli's Theorem~\cite{micchelli1984interpolation} guarantees that the interpolation matrix in~\autoref{augmentedmatrix} is non-singular for various classes of RBFs when distinct nodes are used in the dataset. While non-singularity ensures a unique solution, it does not address the matrix conditioning or the accuracy of the interpolant. Therefore, the selection of $\varepsilon$ remains critical in the sense of balance between accuracy and numerical stability for infinitely smooth RBFs such as GA and MQ.

\subsection{Modified Hermite Radial Basis Functions}
The MHRBF interpolant attempts to enhance interpolation accuracy by incorporating monomial scaling terms into the kernel. This modification mitigates the sensitivity to $\varepsilon$ and reduces the need for higher-order derivatives of the kernel. Specifically the interpolant is now defined as:
\begin{equation}\label{xx}
    s(\vec{x})=\sum^{N}_{i=1}\left\{w_{i}\left[\prod_{j=1}^{d} \left( \vec{x}_{j}-\vec{x}_{i,j} \right)^n\right]+(\vec{x}-\vec{x}_i)^{2n}\cdot\vec{b}_i\right\}\phi(\|\vec{x}-\vec{x}_i\|) + \sum^{M}_{k=1}\lambda_{k}p_k(\vec{x}),
\end{equation}
where $w_i$, $\vec{b}_i$, and $\lambda_i$ have the same meaning as in~\autoref{hermite-first}, powers are performed component-wise, and $d$ represents the dimension of the domain. For clarity, in two-dimensions, the first portion of~\autoref{xx} can be written as
\begin{equation}
    \sum_{i=1}^N\left((x-x_i)^n(y-y_i)^nw_i + (x-x_i)^{2n}\alpha_i + (y-y_i)^{2n}\beta_i\right)\phi(\|\vec{x}-\vec{x}_i\|).
\end{equation}
Using the same condition on the polynomial augmentation, \autoref{eq:polyAugConstraint}, results in the following linear system:
\begin{equation}\label{augmentedmatrix2}
\begin{bmatrix}
    \vec{M}\vec{A} & \vec{N}\vec{A} & \vec{Q}\vec{A} & \vec{P}\\
    \vec{M}_{x}\vec{A} + \vec{M}\vec{A}_{x} & \vec{N}_{x}\vec{A} + \vec{N}\vec{A}_{x} & \vec{Q}\vec{A}_{x} & \vec{P}_x\\
    \vec{M}_{y}\vec{A} + \vec{M}\vec{A}_{y} & \vec{N}\vec{A}_{y} & \vec{Q}_{y}\vec{A} + \vec{Q}\vec{A}_{y} & \vec{P}_y\\
    \vec{P}^T & \vec{P}^T_x & \vec{P}^T_y & \vec{0}
\end{bmatrix}
\begin{bmatrix}
   \vec{w} \\
   \vec{\alpha}\\
   \vec{\beta}\\
   \vec{\lambda}
\end{bmatrix}
=
\begin{bmatrix}
    \vec{f}\\
    \vec{f}_x\\
    \vec{f}_y\\
    \vec{0}
\end{bmatrix},
\end{equation}
where \( \vec{M} \), \( \vec{N} \), and \( \vec{Q} \) are matrices containing the terms \( (x - x_i)^n (y - y_i)^n \), \( (x - x_i)^{2n} \), and \( (y - y_i)^{2n} \), respectively. Although the system is no longer symmetric, the MHRBF only requires first derivatives of the kernel which scale as $h^2$ if $\varepsilon\propto h$, where $h$ is the characteristic spacing between nodes. The HRBF, on the other hand, requires second derivatives that scale as $h^4$. This modification will improve the stability of the method as $h$ decreases. Moving forward, we focus on the GA kernel and leave explorations of other kernels to future work.

\section{Examining the role of Monomial Scaling Term on MHRBF Accuracy}
To facilitate interpretation and focus on the core effects of the monomial scaling term on the behavior of the RBF kernel, we consider the 1D case for the GA kernel. The standard form of the 1D GA function is given by:
\begin{equation}\label{GA1d}
\phi(x)=e^{-(\varepsilon(x-x_0))^2}
\end{equation}
where $\varepsilon$ is the shape parameter controlling the spread of function.
In the MHRBF framework the GA kernel is modified by introducing a monomial term $(x-x_0)^n$, resulting in the function:
\begin{equation}\label{MGA1d}
\phi^*(x)=(x-x_0)^n\phi(x)
\end{equation}
This additional term alters the distribution of the GA kernel by amplifying or suppressing values based on the degree of monomial $n$. Depending on the choice of $n$, this can shift the peak, introduce asymmetry, or modify the decay rate. To better understand the influence of $n$, we now analyze how the degree of the monomial term affects the behavior of the GA kernel.

\subsection{Behavior of the Modified Gaussian Kernel with Varying Degrees of the Monomial \texorpdfstring{$n$}{n}}\label{optimaln1}

The effect of the monomial term is most clearly observed through the choice of $n$. For simplicity, the analysis is restricted to even values of $n$, which eliminates the asymmetry introduced by odd degrees. Our investigations suggest that the odd values of $n$ do not provide significant additional insights beyond those observed with the even values.

To visualize the impact of different monomial degrees, we examine how $\phi^*(x)$ behaves for various values of $n$. We begin by plotting the standard GA kernel and then illustrate how multiplying it by $(x-x_0)^n$ alters the shape and spread of the function. We focus on $\varepsilon=1$ as a mid-range value of the shape parameter to isolate the influence of $n$ without introducing additional complexities. A more detailed analysis of the effects of small and large shape parameter values will be presented in \autoref{optimaln2}. The following plots provide a clear understanding of how the monomial term influences the GA kernel.
\begin{figure}[H]
    \centering
    \begin{subfigure}[t]{0.32\textwidth}
        \centering
        \includegraphics[width=\textwidth]{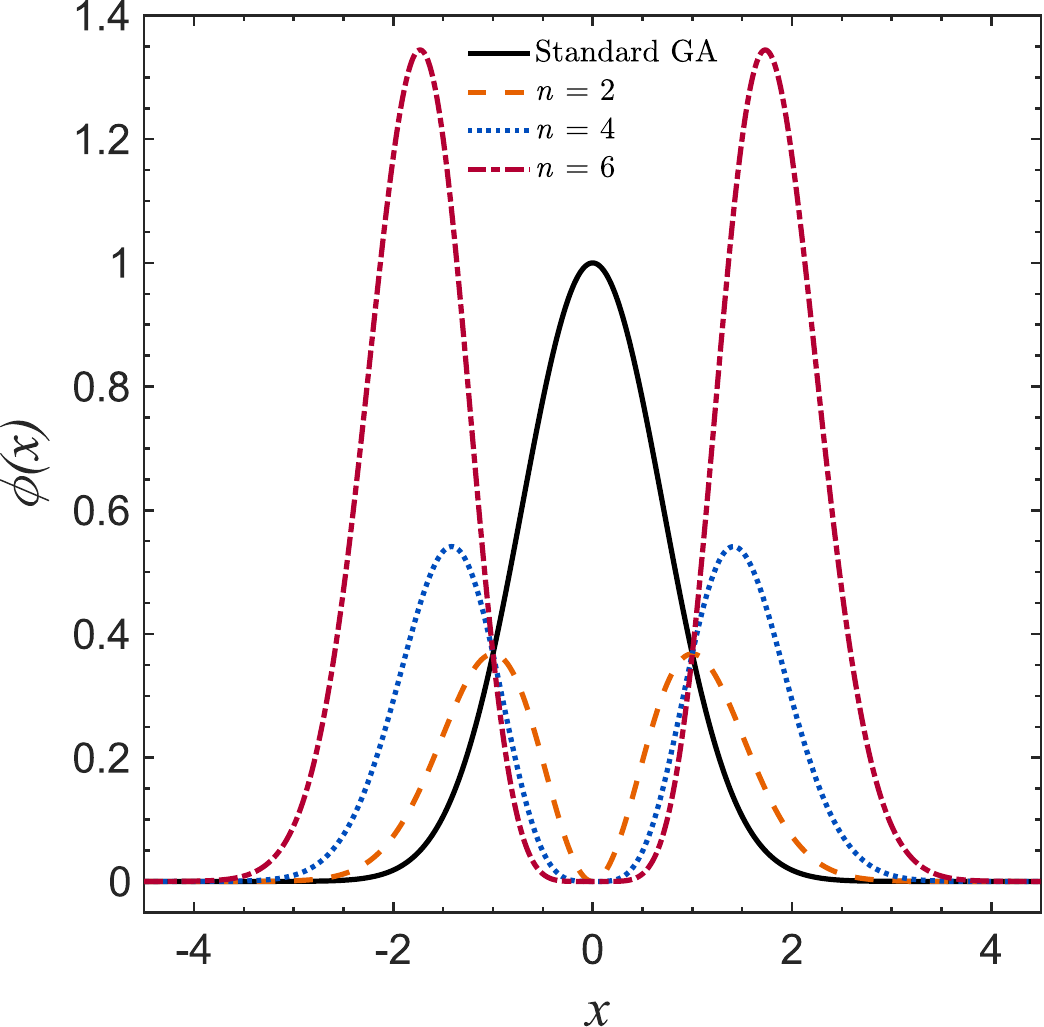}
        \caption{Low-to-mid values of $n$}
        \label{GAleft}
    \end{subfigure}
    \hspace{0.005\textwidth}
    \begin{subfigure}[t]{0.32\textwidth}
        \centering
        \includegraphics[width=\textwidth]{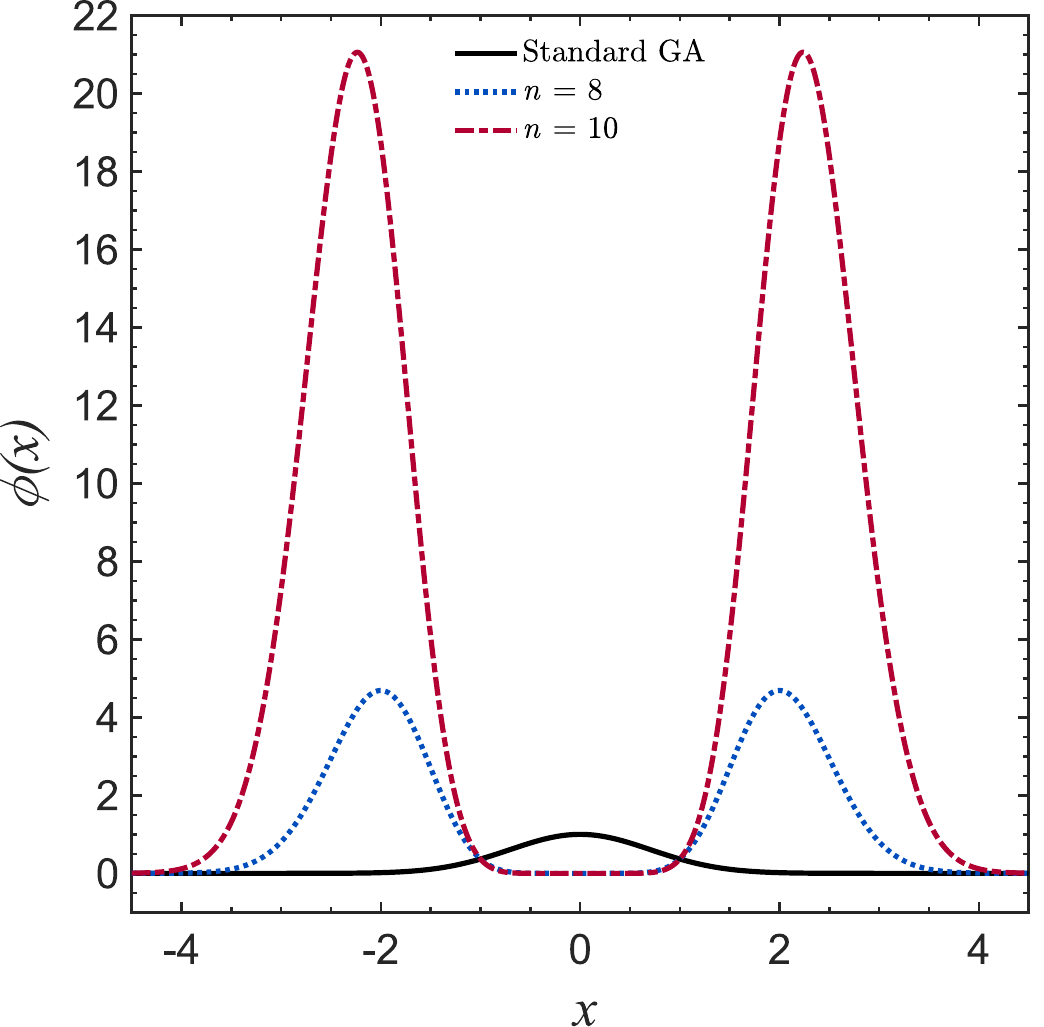}
        \caption{Higher values of $n$}
        \label{GAright}
    \end{subfigure}
    \captionsetup{font=small}
    \caption{Effect of the monomial term on the GA kernel for different values of $n$ using $x_0=0$.}
    \label{GA}
\end{figure}

In \autoref{GA}, we analyze the effect of the monomial scaling term by considering two cases: low-to-mid values of $n$ (\autoref{GAleft}) and higher values of $n$ (\autoref{GAright}). Our observations indicate that the monomial term dictates the location of the kernel's peak and its degree of localization or spread. For a mid-range shape parameter ($\varepsilon\approx1$), increasing $n$ shifts the peak outward and broadens the kernel's effective radius, thereby incorporating more distant nodes. Furthermore, higher values of $n$ lead to an increase in the peak amplitude.

To assess the impact of these modifications, we must determine an optimal $n$ that balances the localization and coverage of the kernel. A highly localized kernel (e.g., $n=2$) concentrates its influence on nearby points (\autoref{GAleft}), which can be beneficial for capturing fine-scale variations. However, this may lead to underfitting if the underlying function exhibits broader trends. Conversely, a broader kernel (e.g., $n=8,10$) can capture those broader trends by considering more distant points (\autoref{GAright}), which might prevent underfitting but could risk overfitting if there is noise in the data.

A notable consequence of large $n$ is a sharp increase in the kernel's peak amplitude. A higher peak amplifies the influence of certain interpolation points, leading to excessive localization and reduced weighting of other data points. This behavior is evident in \autoref{GAright}, where the peak reaches 4.68 for $n=8$ and 21.05 for $n=10$. In this case, the modified kernel behaves similarly to a standard GA kernel with an extremely large shape parameter, resulting in highly localized interpolation near mid-range points and missing interpolation data near center point.

Based on our analysis, mid-range values of $n$ (e.g., $n=4,6$) offer a favorable trade-off between local detail preservation and broader function approximation, essential for stable interpolation. These values maintain sufficient localization while preserving an appropriate kernel radius, making them a practical choice for robust interpolation.

\subsection{Exploring the Optimal Combination of \texorpdfstring{$\varepsilon$}{e} and \texorpdfstring{$n$}{n} for MHRBF Accuracy}
\label{optimaln2}

\begin{figure}[H]
    \centering
        \begin{subfigure}{0.24\textwidth}
            \centering
            \includegraphics[width=\linewidth]{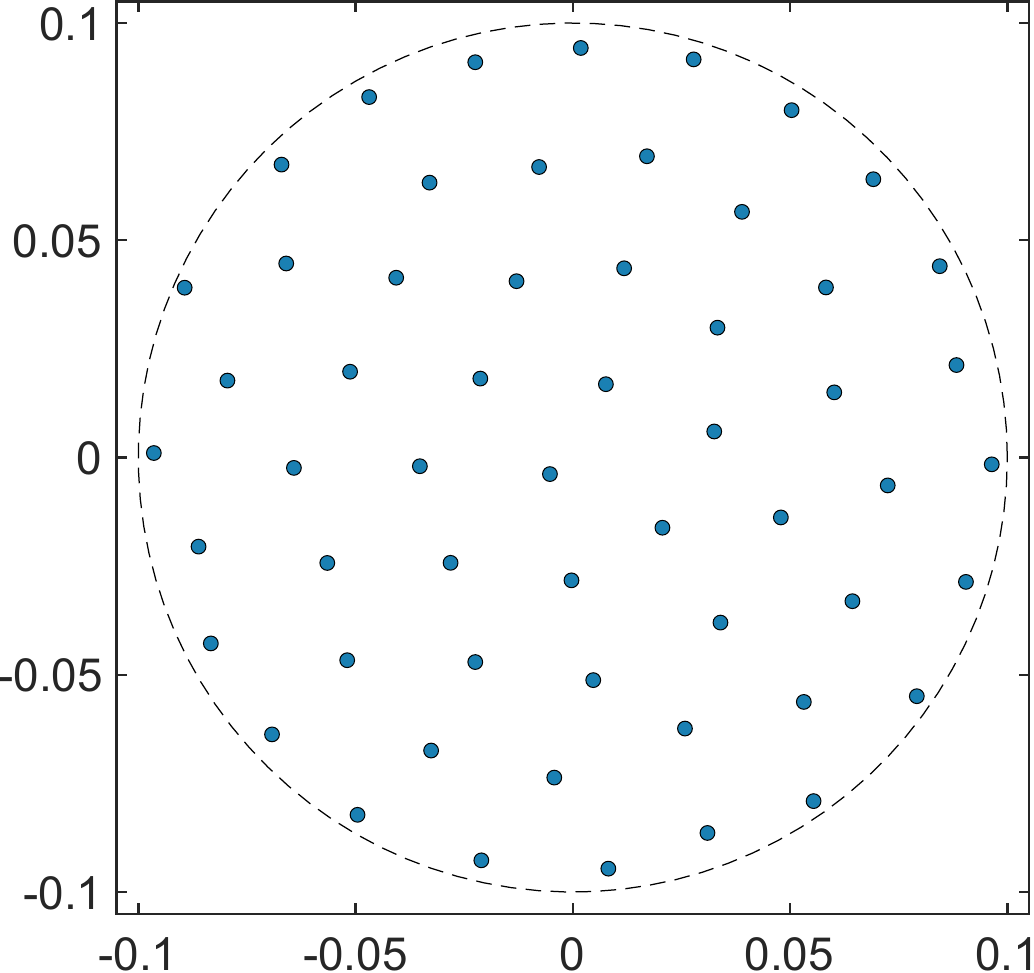}
            \captionsetup{justification=centering}
            \caption{Data nodes}
        \end{subfigure}
        \hspace{0.01\textwidth}
        \begin{subfigure}{0.24\textwidth}
            \centering
            \includegraphics[width=\linewidth]{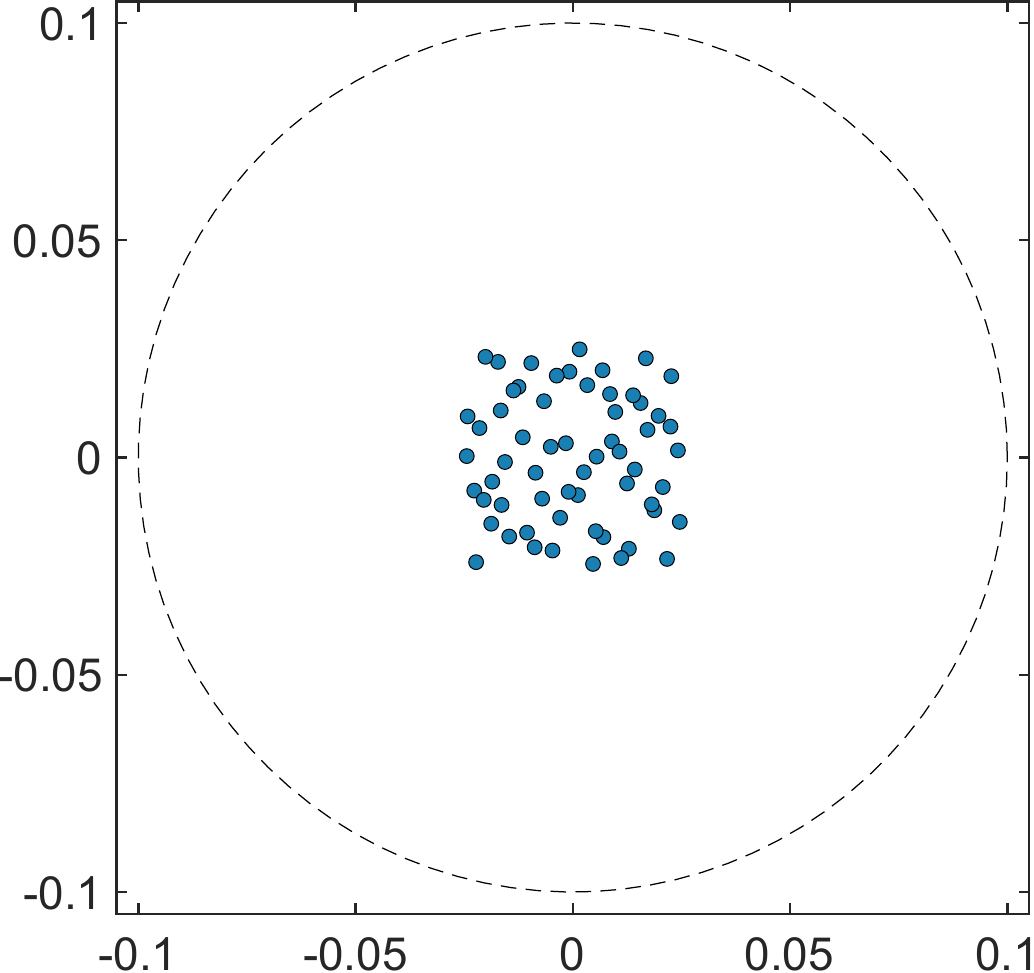}
            \captionsetup{justification=centering}
            \caption{Error evaluation nodes}
        \end{subfigure}
    \captionsetup{font=small}
    \caption{Distribution of data and error evaluation nodes within a circle of radius 0.1.}
    \label{minenergy}
\end{figure}

To evaluate the interaction of $\varepsilon$ and $n$ on interpolation accuracy of the MHRBF method, we implement the method in Matlab using double-precision and utilize a benchmark case inspired by Ref.~\cite{flyer2016role} using the GA RBF kernel. The data given consists of 56 nodes arranged in a minimum energy configuration~\cite{Moreno2025} contained in a radius of $R=0.1$. The errors are computed using 60 Halton set evaluation nodes contained near the center, see \autoref{minenergy}. In all cases standard double-precision is used.

The function used throughout this analysis is defined as:
\begin{equation}\label{function1}
f(x,y)=\sin(6x)+\cos(4y)+\sin(3x+2y).
\end{equation}
To assess the accuracy of the MHRBF method, we compute the $L_\infty$ error for function interpolation and its first derivatives in both the $x$ and $y$ directions over a range of $\varepsilon$ and $n$ values. Specifically, we consider $\varepsilon \in[0.001,10]$ and $n \in[1,9]$. The results are visualized using contour plots, where the $x$-axis represents the shape parameter $\varepsilon$ and the $y$-axis corresponds to the monomial degree $n$. These plots reveal optimal parameter combinations that minimize interpolation error.

\begin{figure}[H]
    \centering
    \begin{subfigure}{0.32\textwidth}
        \centering
        \includegraphics[width=\linewidth]{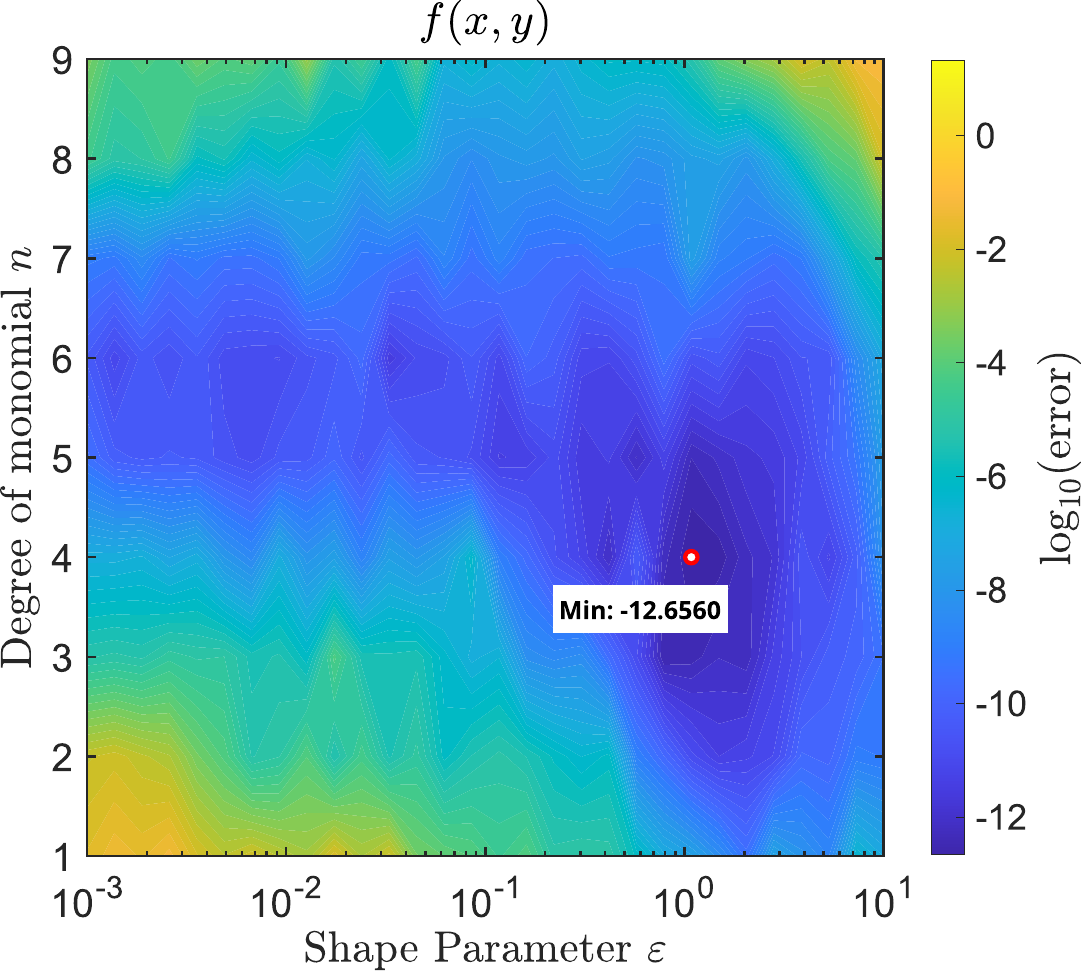}
        \caption{$\log_{10}(\text{error})$ for $f$.}
        \label{contourf}
    \end{subfigure}
    \hspace{0.005\textwidth}
    \begin{subfigure}{0.32\textwidth}
        \centering
        \includegraphics[width=\linewidth]{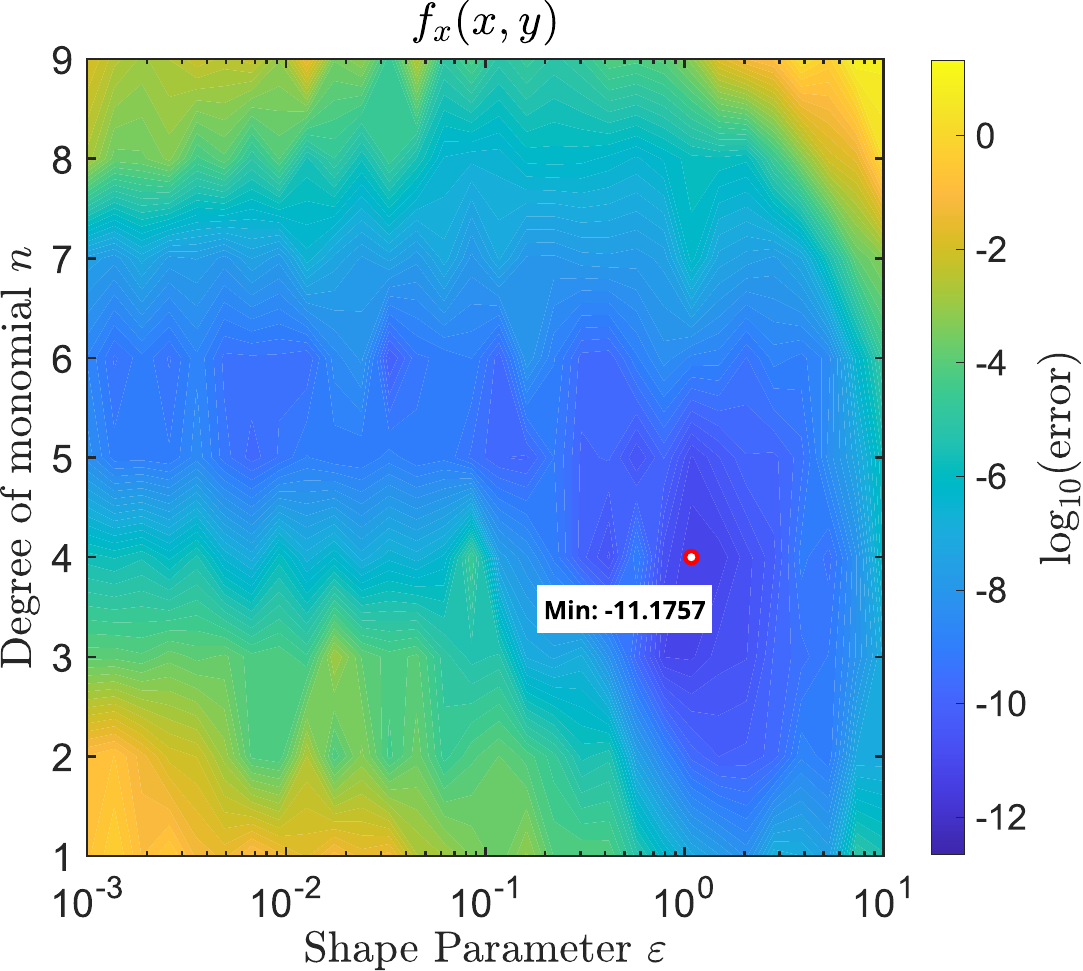}
        \caption{$\log_{10}(\text{error})$ for $f_x$.}
        \label{contourfx}
    \end{subfigure}
    \hspace{0.005\textwidth}
    \begin{subfigure}{0.32\textwidth}
        \centering
        \includegraphics[width=\linewidth]{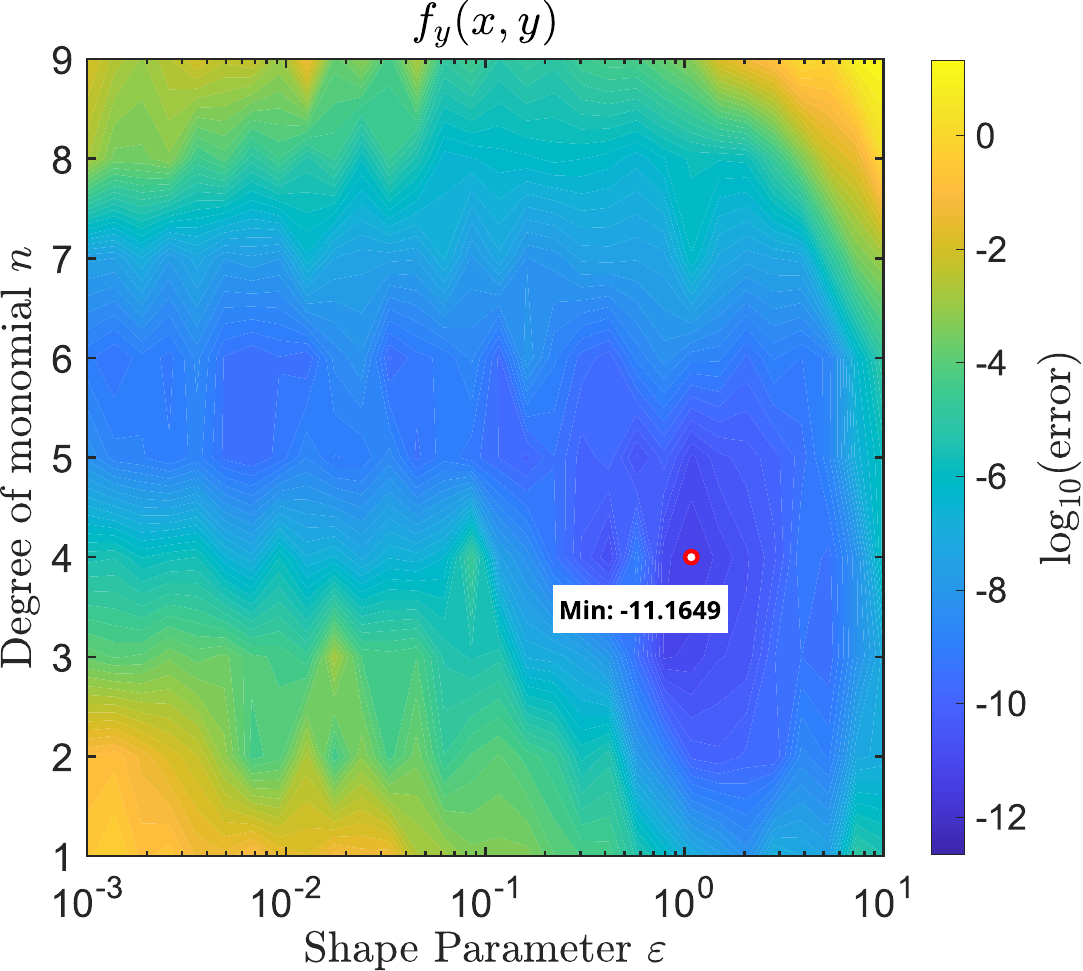}
        \caption{$\log_{10}(\text{error})$ for $f_y$.}
        \label{contourfy}
    \end{subfigure}

    \captionsetup{font=small}
    \caption{Contour plots of the $L_\infty$ error for the MHRBF method as a function of shape parameter $\varepsilon$ and monomial degree $n$. Separate plots are shown for function values and first derivatives ($f_x$, $f_y$), revealing optimal parameter combinations that minimize error.}
    \label{Contour}
\end{figure}

The contour plots in~\autoref{Contour} illustrate how the interpolation error varies across different combinations of $\varepsilon$ and $n$. For small shape parameters ($\varepsilon\ll1$), the Gaussian kernel exhibits a nearly uniform spread, leading to a slow decay across the field. In this regime, higher degree monomial terms grow rapidly, which helps counteract the kernel's flatness by introducing localized weighting effects. This adjustment enhances the kernel's ability to capture fine scale variations in the function. As observed in~\autoref{Contour} for $\varepsilon = 10^{-3}$ the optimal monomial degree that minimizes the interpolation error is $n=6$. However, as discussed in \autoref{optimaln1}, choosing $n$ values that are either too small or too large can disrupt the balance between localization and coverage, leading to a loss of accuracy. This effect is evident in the contour plots, where deviations from $n=6$ in either direction result in an increase in error.

For mid-range shape parameters ($\varepsilon \approx 1$), the standard Gaussian RBF kernel exhibits a more reasonable decay, reducing the need for high degree monomial scaling. Consequently, the optimal $n$ is lower compared to the small $\varepsilon$ case. The results indicate that for $\varepsilon = 1$, the optimal monomial degree is $n=4$, which aligns with our findings in \autoref{optimaln1}.

As the shape parameter increases ($\varepsilon \gg 1$), the standard Gaussian RBF kernel becomes highly localized, leading to rapid decay. In this case, a monomial term with a lower degree, which grows more gradually, helps counteract excessive localization and ensures sufficient coverage of the interpolation domain. This behavior is evident in~\autoref{Contour}, where for large $\varepsilon$, the minimal error is achieved at a lower monomial degree of $n=3$.

These results indicate that choosing a shape parameter around $\varepsilon \approx 1$ and a monomial degree near $n \approx 4$ offers a well-balanced approach, maintaining both localization and coverage for stable and accurate interpolation in the rest of this study.

\section{MHRBF Interpolation with Polynomial Augmentation}\label{polyeffect}
In this section, we investigate the effect of polynomial augmentation on MHRBF interpolation accuracy. As with the Standard RBF and HRBF methods, we expect that the MHRBF method can also benefit from polynomial augmentation, potentially improving accuracy by capturing low-order polynomial trends in the data. Specifically, we examine two cases: augmentation with a linear polynomial (degree 1) and a higher-order polynomial (degree 9). By comparing these cases, we aim to determine how the inclusion of polynomials influences the interpolation error. The analysis follows the same setting as in \autoref{optimaln2}, where we used a GA kernel and interpolated a function sampled at 56 nodes arranged in a minimum energy configuration, with errors evaluated using 60 Halton set points near the center.

To assess the impact of polynomial augmentation, we compute the $L_\infty$ error for function interpolation and its first derivatives in the $x-$ and $y-$directions using double-precision mathematics. As concluded in \autoref{optimaln2}, we select a monomial degree of $n=4$ for this section. The results are presented as error plots, where the $x$-axis represents the shape parameter $\varepsilon$ and the $y$-axis corresponds to the $L_\infty$ error. Each plot contains two curves, corresponding to polynomial degrees 1 and 9, allowing for a direct comparison of their influence on interpolation accuracy.

\begin{figure}[H]
    \centering
    \begin{subfigure}{0.32\textwidth}
        \includegraphics[width=\linewidth]{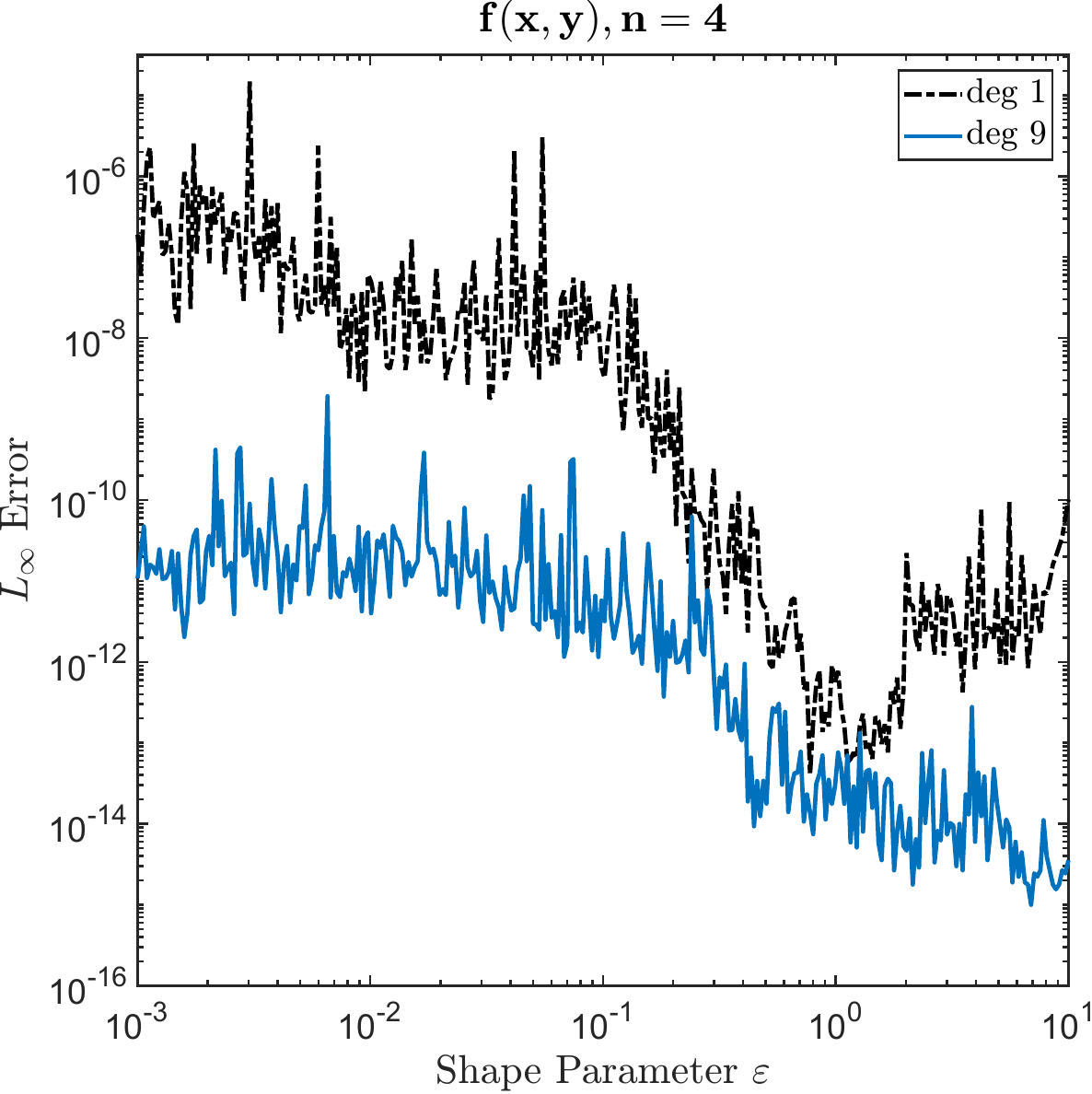}
        \caption{Error in $f$ vs. shape parameter}
    \end{subfigure}
    \hspace{0.005\textwidth}
    \begin{subfigure}{0.32\textwidth}
        \centering
        \includegraphics[width=\linewidth]{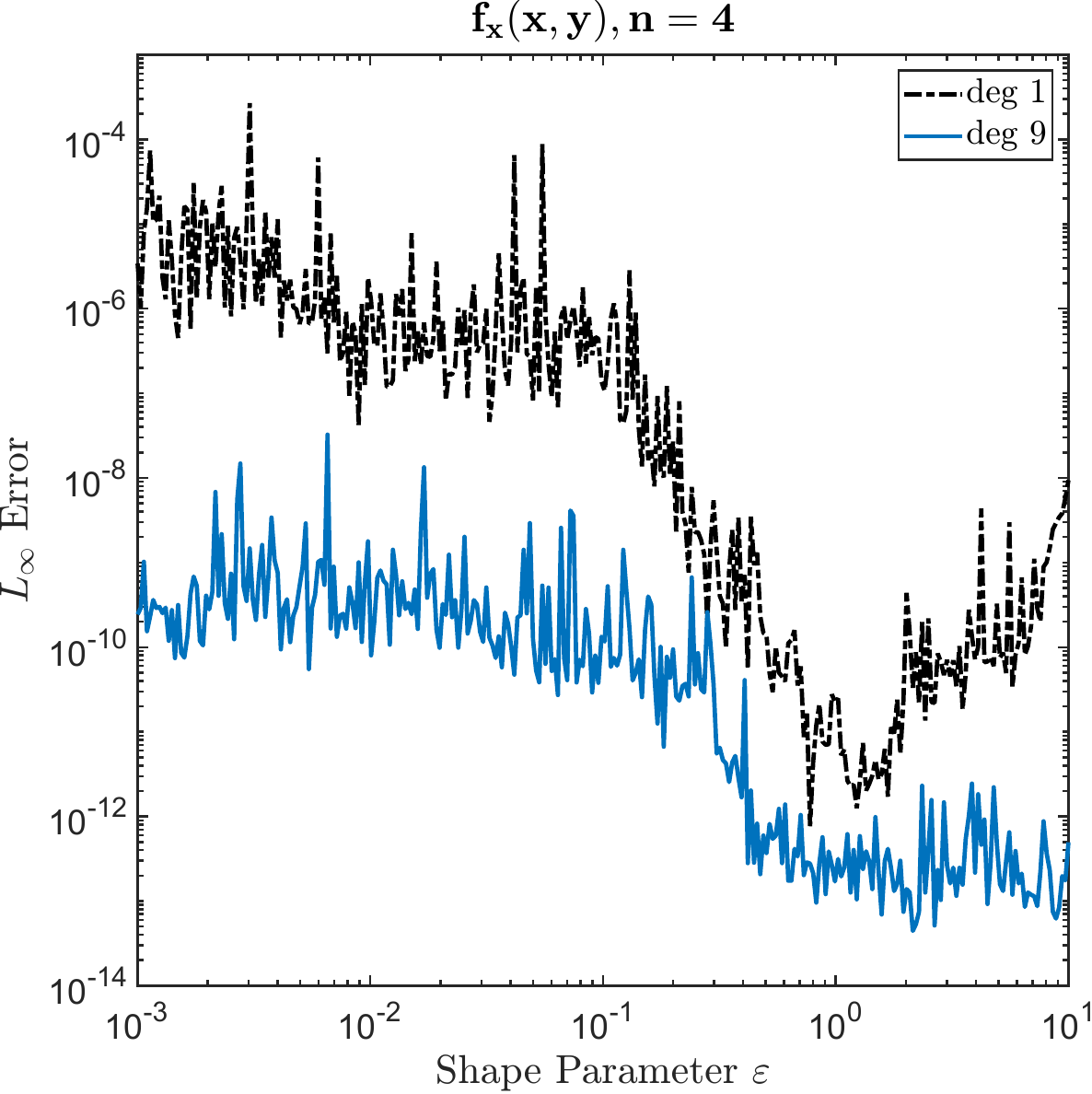}
        \captionsetup{justification=centering}
        \caption{Error in $f_x$ vs. shape parameter}
    \end{subfigure}
    \hspace{0.005\textwidth}
    \begin{subfigure}{0.32\textwidth}
        \centering
        \includegraphics[width=\linewidth]{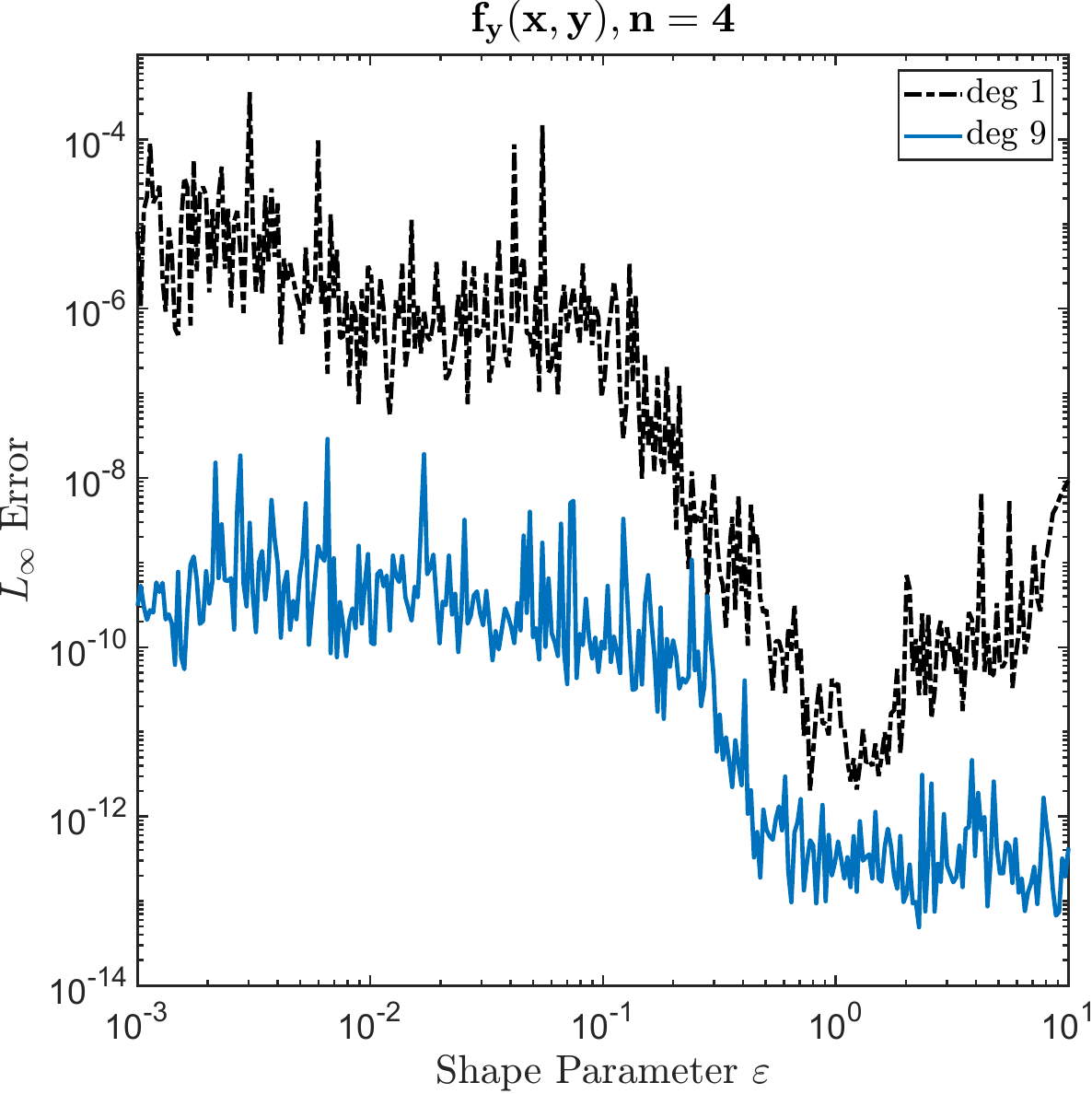}
        \captionsetup{justification=centering}
        \caption{Error in $f_y$ vs. shape parameter}
    \end{subfigure}
    \captionsetup{font=small}
    \caption{Comparison of MHRBF accuracy with polynomial augmentation of degrees 1 and 9: The errors (function, $f$, and its first derivatives, $f_x$ and $f_y$) computed in double precision using the GA kernel, as functions of shape parameter $\varepsilon$, with $n=4$.}
    \label{fig:polyeffect}
\end{figure}
The results of the error plots in \autoref{fig:polyeffect} demonstrate that polynomial augmentation significantly influences the interpolation accuracy of the MHRBF method across different values of the shape parameter $\varepsilon$. In particular, augmentation with a higher degree polynomial (degree 9) consistently reduces the $L_\infty$ error more than a linear polynomial (degree 1), especially at the extremes of the shape parameter range. For example, for function interpolation at $\varepsilon = 0.001$, the error drops from the order of $10^{-7}$ (degree 1) to $10^{-11}$ (degree 9), and at $\varepsilon = 10$, it improves from $10^{-10}$ to $10^{-15}$. A similar trend is observed for the interpolation errors of the first derivatives, $f_x$ and $f_y$.

However, the difference in accuracy between the two polynomial degrees becomes much less significant around $\varepsilon = 1$, where both configurations yield extremely low errors (on the order of $10^{-13}$ for degree 1 and $10^{-14}$ for degree 9). As concluded in \autoref{optimaln2}, this mid-range value of $\varepsilon$ corresponds to the optimal balance between kernel localization and global coverage, resulting in minimal interpolation error.

Given that $\varepsilon \approx 1$ provides near optimal accuracy and the added benefit from increasing the polynomial degree from 1 to 9 is marginal in this region, using a higher degree polynomial is not justified when considering computational cost. Therefore, a linear polynomial (degree 1) offers a more practical and efficient choice for polynomial augmentation in MHRBF interpolation, especially when targeting shape parameters near their optimal range.

\section{HRBF vs. MHRBF}
In this section, we compare the performance of the standard HRBF method with the proposed MHRBF method. The comparison is organized into three parts: interpolation accuracy, the effect of node spacing, and computational cost. This structure allows us to evaluate not only the improvements in approximation quality achieved by the MHRBF method, but also the trade-offs in scalability and efficiency. As before, we utilize double-precision for these calculations.

\subsection{Assessment of Interpolation Accuracy}\label{accuracy}
The first scenario explores the interpolation accuracy of both methods by computing the $L_\infty$ error for function interpolation and its first derivatives in both the $x-$ and $y-$directions. The data and evaluation nodes are the same as \autoref{optimaln2}. To provide a more comprehensive evaluation, we consider two distinct target functions: the first is the function defined in \autoref{function1}, and the second is the Six-Hump Camelback function. The goal is to assess the accuracy benefits of MHRBF across a range of shape parameters ($\varepsilon$), varying from $0.001$ to $10$. Based on the findings in \autoref{optimaln2} and the polynomial augmentation study in \autoref{polyeffect}, we fix the monomial degree at $n = 4$ and use a linear polynomial (degree 1) for augmentation in both methods. We begin by presenting the results for the function in \autoref{function1}, which include plots of the $L_\infty$ error as a function of $\varepsilon$ for both HRBF and MHRBF methods.
\begin{figure}[H]
    \centering
    \begin{subfigure}{0.32\textwidth}
        \includegraphics[width=\linewidth]{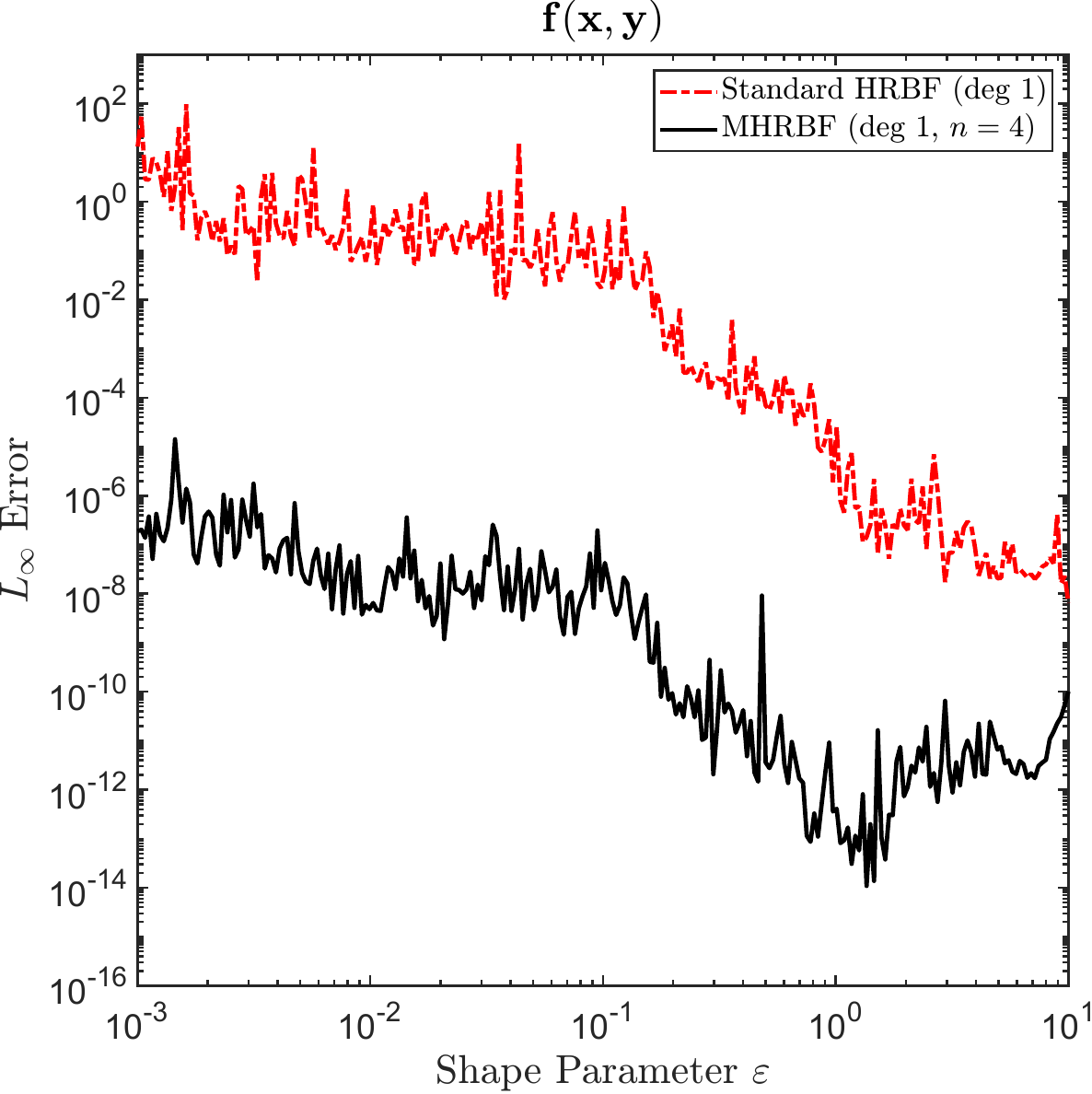}
        \caption{Error in $f$ vs. shape parameter}
        \label{fig:accuracyf1}
    \end{subfigure}
    \hspace{0.005\textwidth}
    \begin{subfigure}{0.32\textwidth}
        \includegraphics[width=\linewidth]{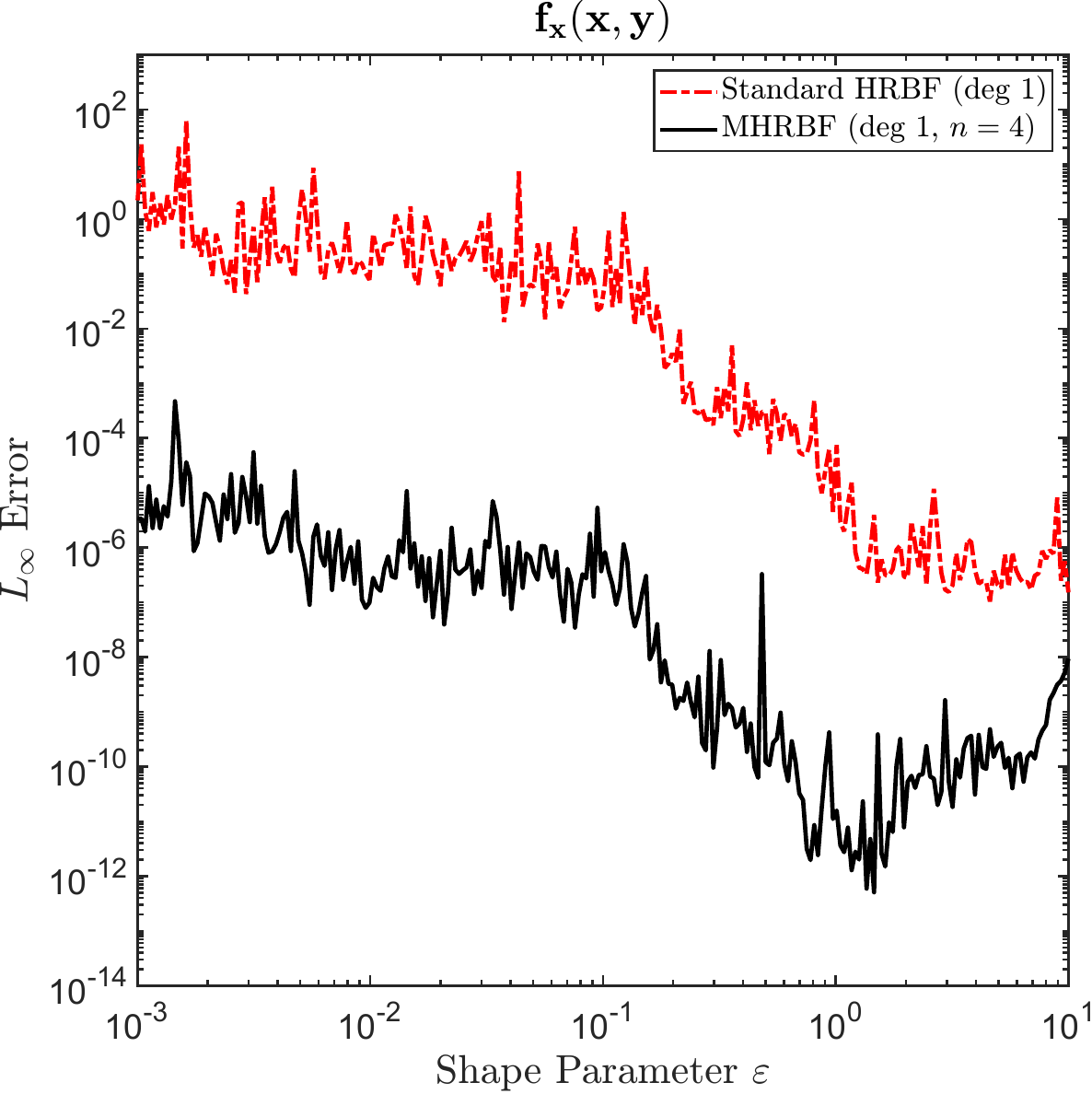}
        \caption{Error in $f_x$ vs. shape parameter}
        \label{fig:accuracyfx1}
    \end{subfigure}
    \hspace{0.005\textwidth}
    \begin{subfigure}{0.32\textwidth}
        \includegraphics[width=\linewidth]{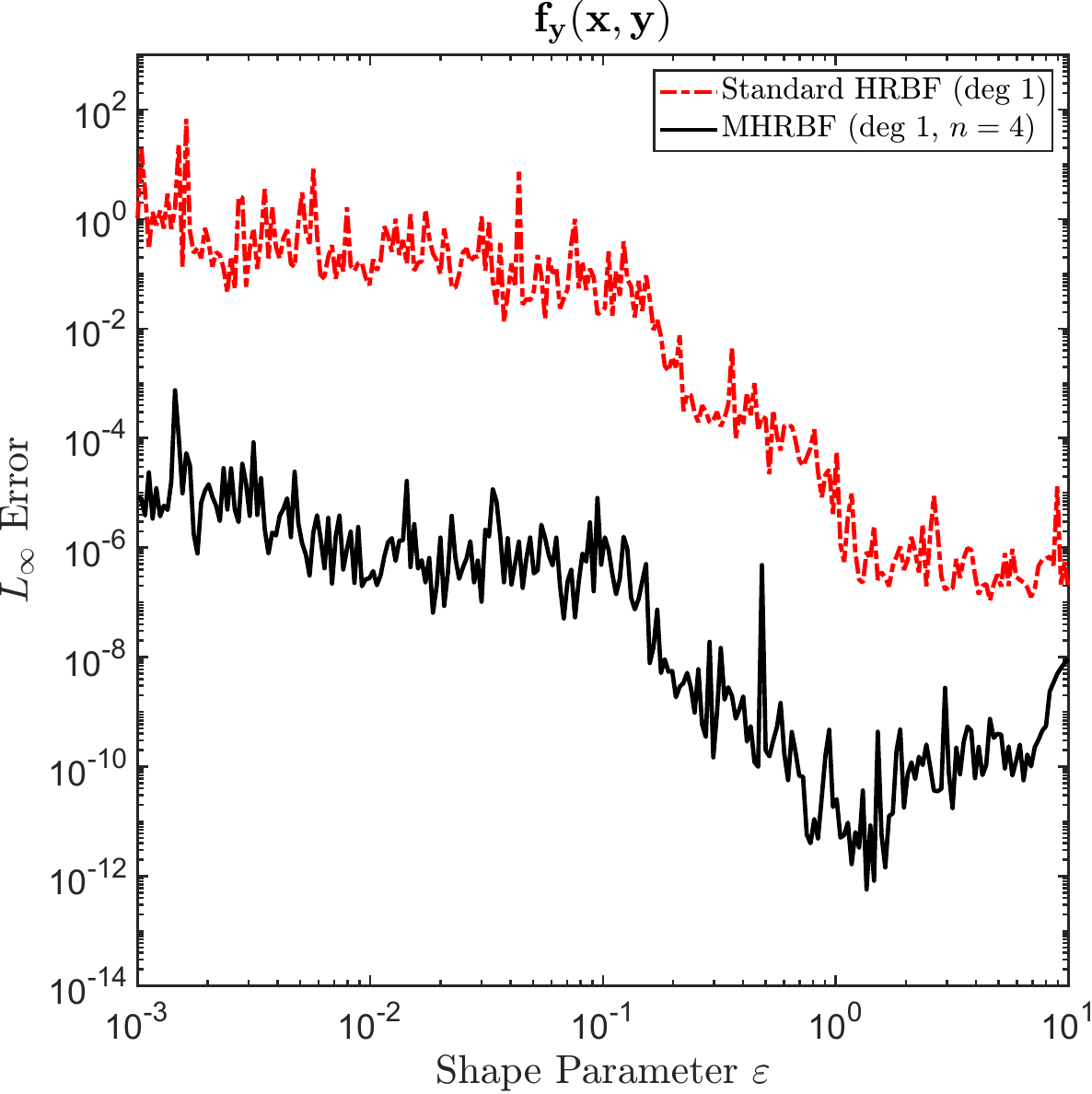}
        \caption{Error in $f_y$ vs. shape parameter}
        \label{fig:accuracyfy1}
    \end{subfigure}
    \captionsetup{font=footnotesize}
    \caption{Comparison of Accuracy for HRBF and MHRBF: The errors (function, $f$, and its first derivatives, $f_x$ and $f_y$) computed in double precision for GA kernel, as functions of shape parameter $\varepsilon$, with polynomial degree of $1$ and $n=4$ for MHRBF. Results correspond to the trigonometric test function.}
    \label{accuracyf1}
\end{figure}

The results in \autoref{accuracyf1} for the trigonometric test function reveal significant accuracy improvements with MHRBF compared to HRBF. As the error trends for $\nabla f$ closely follow that of $f$, the subsequent analysis focuses primarily on the function interpolation.

\autoref{fig:accuracyf1} demonstrates that the MHRBF achieves lower error over all the $\varepsilon$ values compared with the standard HRBF. In the standard HRBF for very small shape parameters ($\varepsilon\ll 1$), the GA kernel becomes nearly flat, making the entries in the upper left block of the system matrix (\autoref{augmentedmatrix}) almost identical to 1. Meanwhile, derivative related terms approach zero, failing to sufficiently break the near linear dependencies in the matrix. Combined with the symmetry of the system matrix, this near linear dependence leads to severe ill-conditioning. The key innovation of MHRBF lies in its polynomial scaling terms, which introduce asymmetry by weighting points differently based on their position relative to ($x_0,y_0$). Consequently, the MHRBF addresses the issue of symmetry and near linear dependencies in standard HRBF, stabilizing the matrix and improving accuracy for very small $\varepsilon$. For moderate values of the shape parameter ($\varepsilon\approx1$), in addition to the symmetric structure of the system matrix (\autoref{augmentedmatrix}), the dominant kernel terms in the upper left block still overshadow the derivative terms. However, the polynomial scaling terms in MHRBF mitigate this imbalance, allowing derivative terms to play a more proportional role in interpolation, thus preserving accuracy. For large shape parameters ($\varepsilon\approx10$), both methods experience accuracy degradation due to underfitting as the GA kernel and its derivatives rapidly decay to zero, leading to sparse interactions between nodes and reduced interpolation quality.

\begin{figure}[H]
    \centering
    \begin{subfigure}{0.32\textwidth}
        \includegraphics[width=\linewidth]{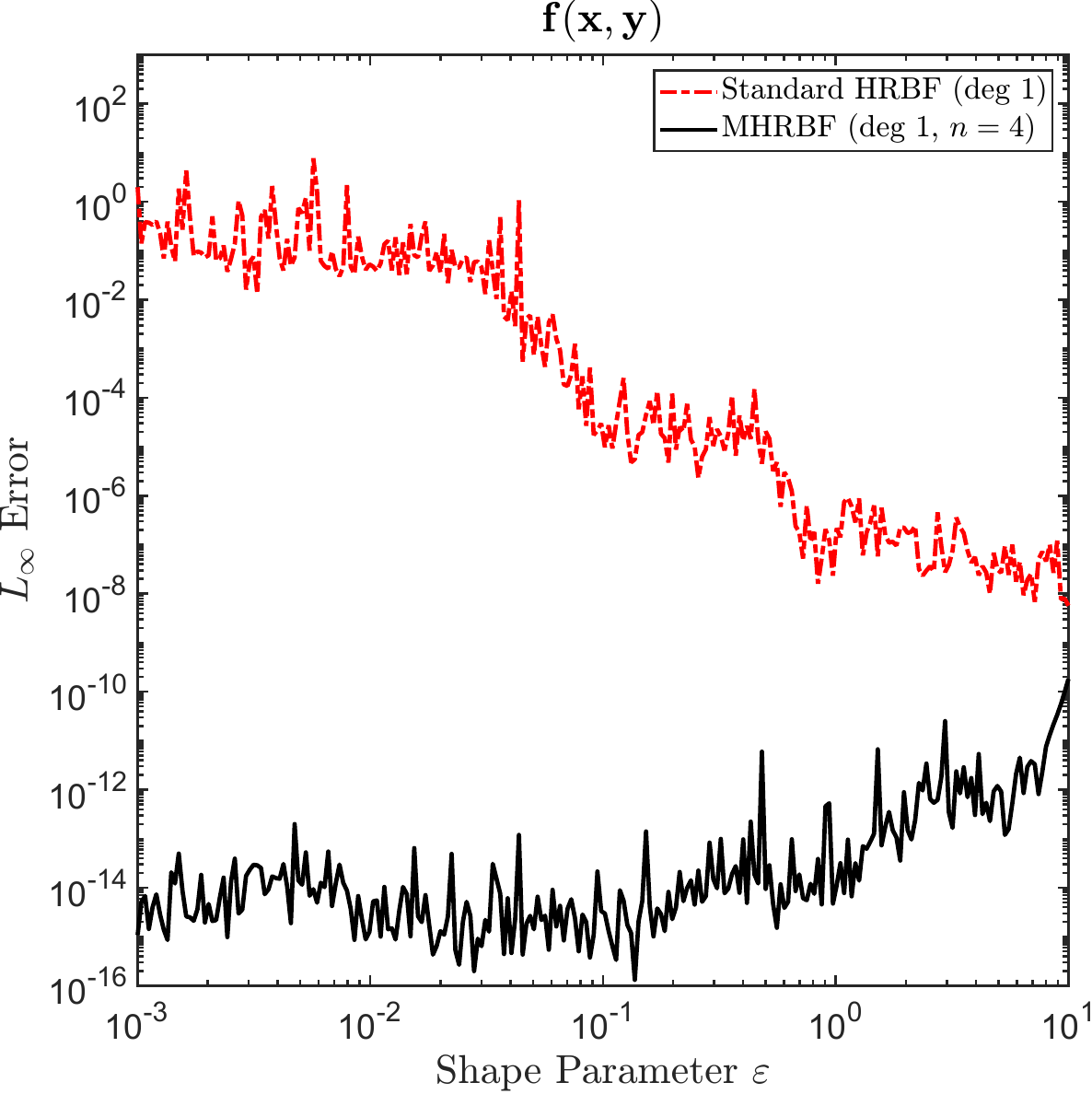}
        \caption{Error in $f$ vs. shape parameter}
        \label{fig:accuracyf2}
    \end{subfigure}
    \hspace{0.005\textwidth}
    \begin{subfigure}{0.32\textwidth}
        \includegraphics[width=\linewidth]{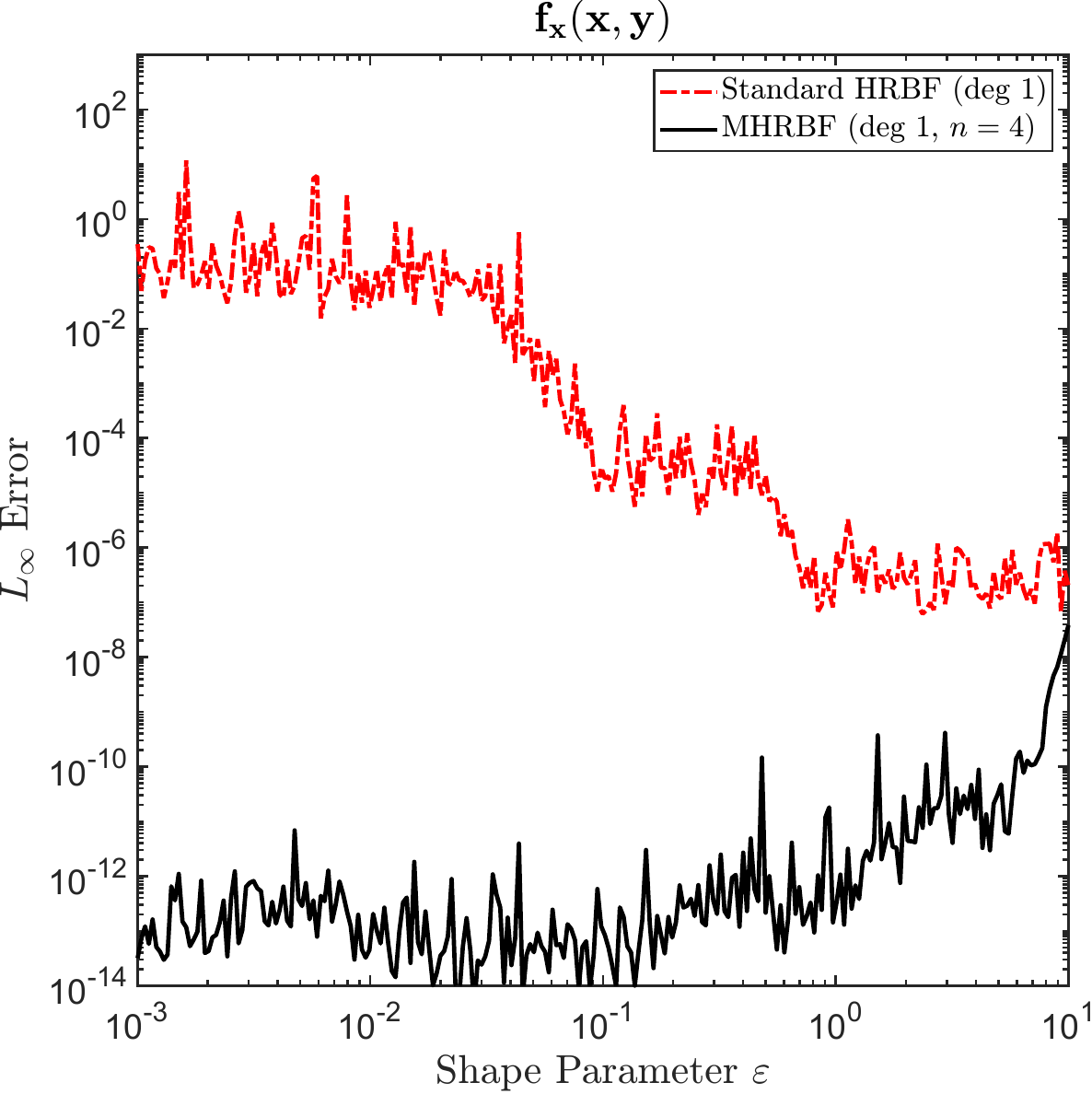}
        \caption{Error in $f_x$ vs. shape parameter}
        \label{fig:accuracyfx2}
    \end{subfigure}
    \hspace{0.005\textwidth}
    \begin{subfigure}{0.32\textwidth}
        \includegraphics[width=\linewidth]{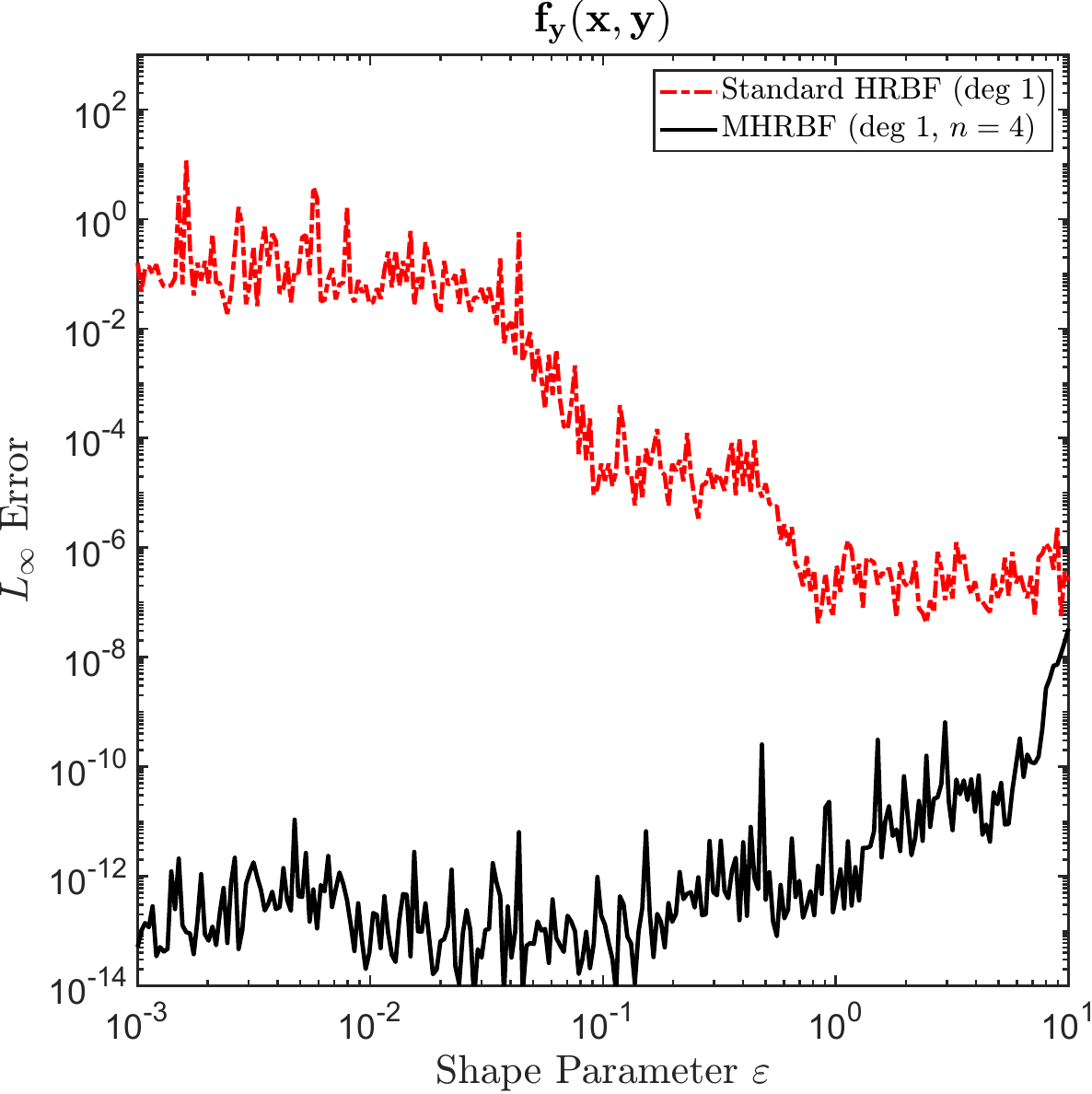}
        \caption{Error in $f_y$ vs. shape parameter}
        \label{fig:accuracyfy2}
    \end{subfigure}

    \captionsetup{font=footnotesize}
    \caption{Comparison of Accuracy and Numerical Stability for HRBF and MHRBF: The errors (function, $f$, and its first derivatives, $f_x$ and $f_y$) and condition number computed in double precision for GA kernel, as functions of shape parameter $\varepsilon$, with polynomial degree of $6$ and $n=4$ for MHRBF. Results correspond to the Six-Hump Camelback function.}
    \label{accuracyf2}
\end{figure}

To further validate the improved accuracy of the MHRBF method over the standard HRBF method, we repeat the interpolation test using a second target function: the well-known Six-Hump Camelback function, defined as:
\begin{equation}\label{six-hump}
f(x, y) = \left(4 - 2.1x^2 + \frac{x^4}{3}\right)x^2 + xy + (-4 + 4y^2)y^2.
\end{equation}
This function introduces polynomial nonlinearity and mixed terms, providing a contrasting test case to the function analyzed earlier. The corresponding results using an augmenting polynomial of order 6 are presented in \autoref{accuracyf2}.

The results in \autoref{accuracyf2} confirm the superior accuracy of the MHRBF method across the entire range of $\varepsilon$ values for the Six-Hump Camelback function as well. Compared to the trigonometric case, the behavior of the interpolation error exhibits a flatter profile for small to moderate values of the shape parameter. Specifically, the MHRBF error for the function (\autoref{accuracyf2}) starts below $10^{-12}$ for $\varepsilon = 10^{-3}$ and remains in the range of $10^{-14}$ to $10^{-11}$ until approximately $\varepsilon = 0.5$, after which it increases gradually, reaching around $10^{-10}$ at $\varepsilon = 10$. In contrast, the HRBF error begins near $10^{-1}$, decreases to about $10^{-7}$, and then fluctuates between $10^{-8}$ and $10^{-6}$ for larger $\varepsilon$ values. Unlike the previous function, where MHRBF error initially decreased and then increased beyond $\varepsilon \approx 1$, the current results reveal a remarkably flat error trend at low-to-mid $\varepsilon$ values, indicating high robustness of MHRBF to shape parameter variation in this case. These findings further reinforce the improved stability and approximation quality of the MHRBF method, particularly for functions with polynomial and mixed nonlinearities.

\subsection{Convergence Behavior with Respect to Node Spacing}\label{convergence}

\begin{figure}[H]
    \centering
    \begin{subfigure}{0.32\textwidth}
        \includegraphics[width=\linewidth]{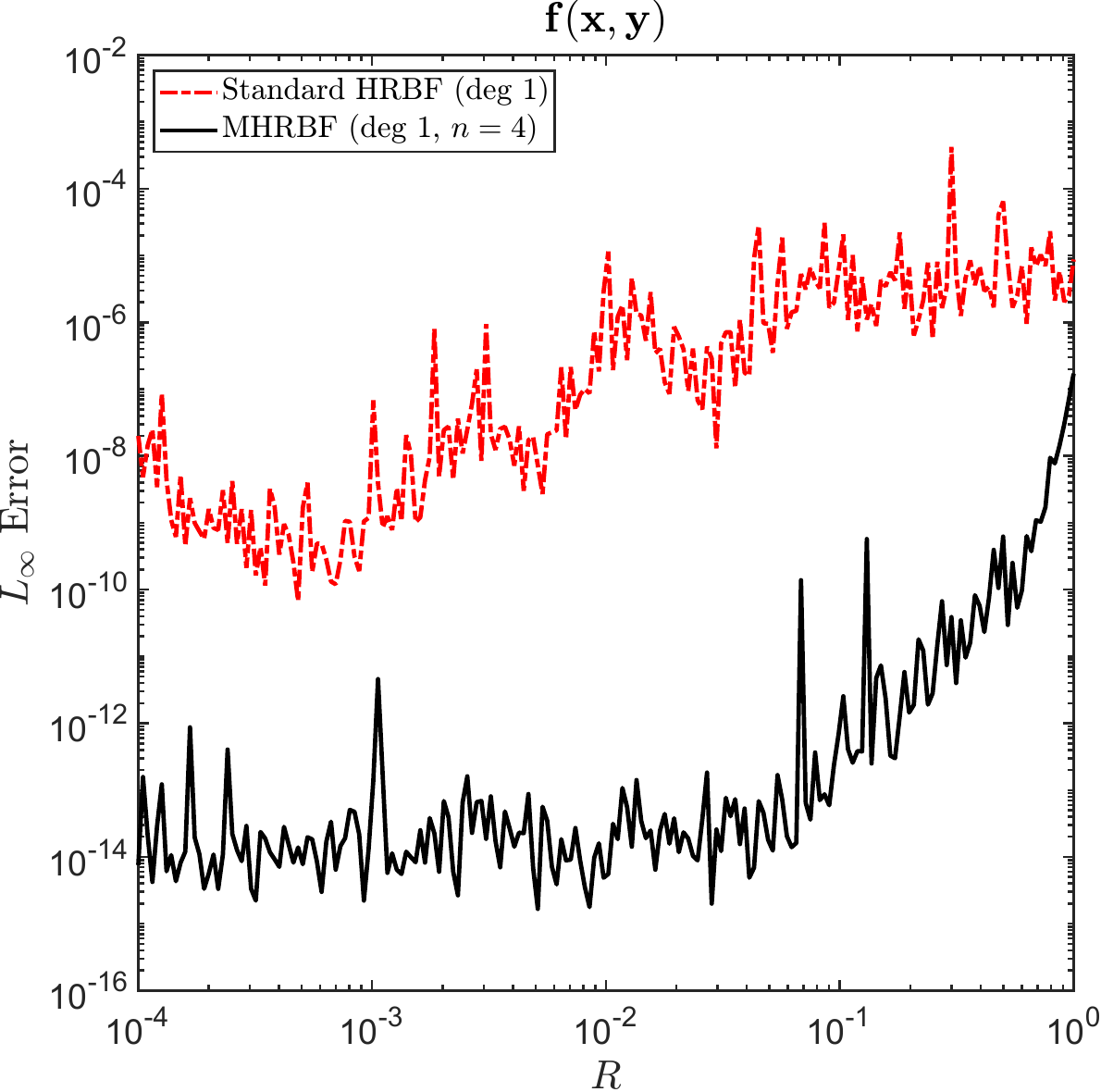}
        \caption{Error in $f$ vs. radius}
        \label{fig:error_f}
    \end{subfigure}
    \hspace{0.005\textwidth}
    \begin{subfigure}{0.32\textwidth}
        \includegraphics[width=\linewidth]{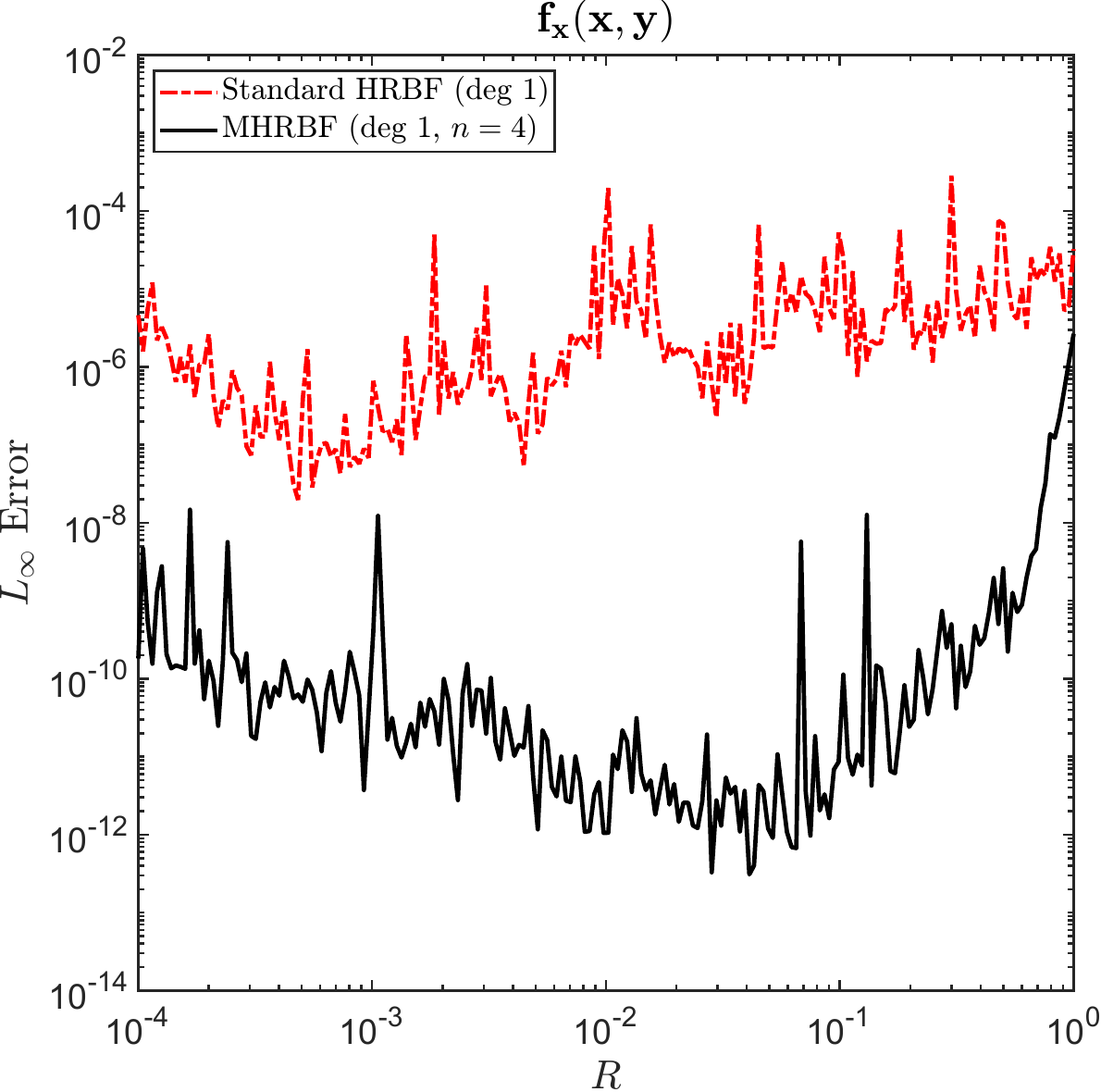}
        \caption{Error in $f_x$ vs. radius}
    \end{subfigure}
    \hspace{0.005\textwidth}
    \begin{subfigure}{0.32\textwidth}
        \includegraphics[width=\linewidth]{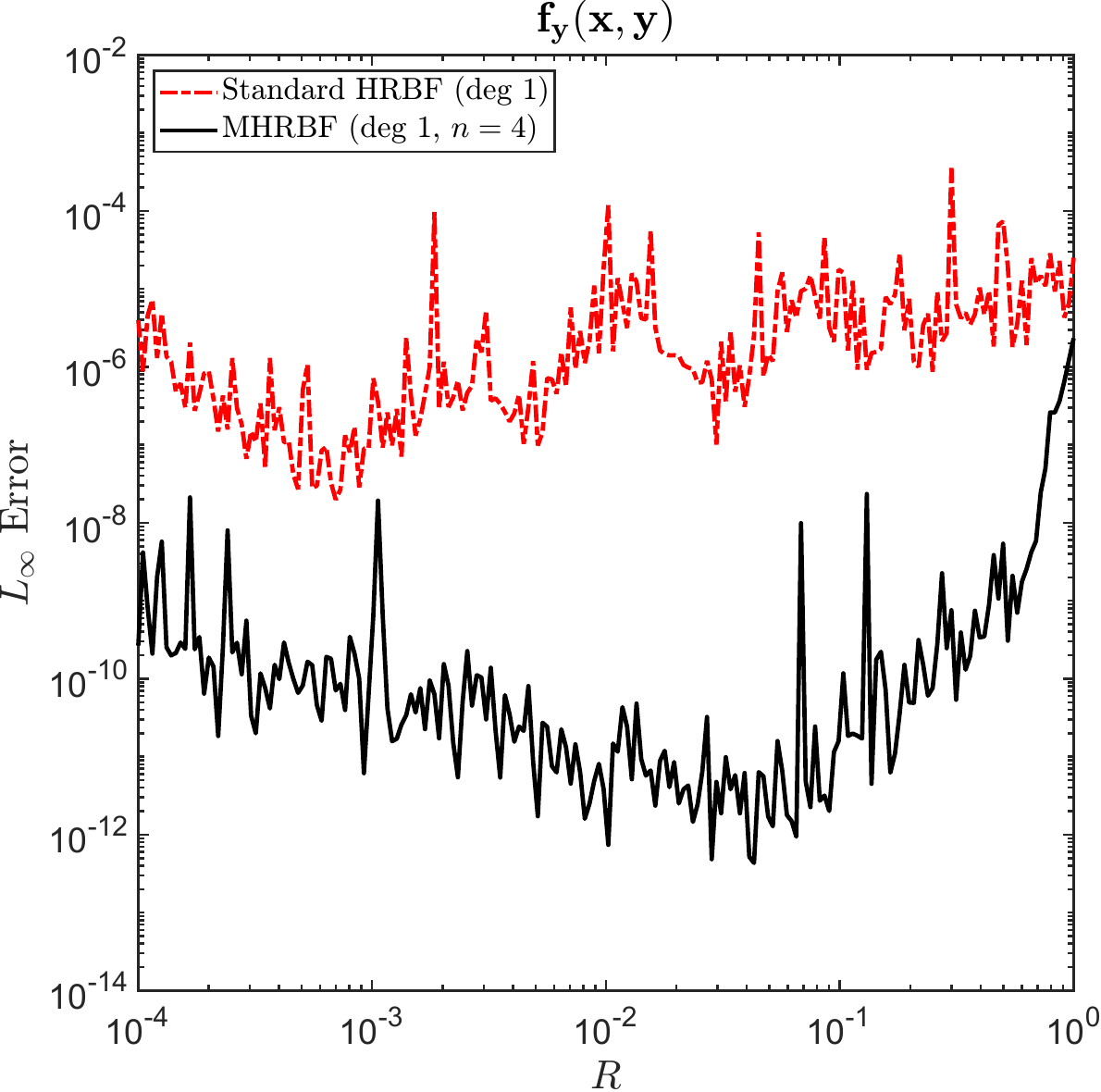}
        \caption{Error in $f_y$ vs. radius}
        \label{fig:cond}
    \end{subfigure}

    \captionsetup{font=footnotesize}
    \caption{The interpolation errors (function, $f$, and its first derivatives, $f_x$ and $f_y$) computed in double precision using the GA kernel, shown as functions of the interpolation radius $R$ (i.e., node spacing). The shape parameter is fixed at $\varepsilon = 1$, with a polynomial degree of $1$ and $n = 4$ for MHRBF. Results correspond to the trigonometric function.}
    \label{scalingResults}
\end{figure}

To investigate the effect of node spacing on interpolation performance, we perform a scaling analysis by varying the radius of the interpolation domain, $R$, over several orders of magnitude from $10^{-4}$ to $10$. In this setting, the relative positions of the data and evaluation nodes are fixed, and the radius acts as a uniform scaling factor for their spatial distribution. This allows us to examine how the distance between the interpolation nodes influences the accuracy and convergence behavior. The shape parameter $\varepsilon$ is kept fixed at 1, the augmenting polynomial order is 1, and the interpolation error is evaluated as a function of $R$. The results are summarized in \autoref{scalingResults}.
For both methods, the convergence plots exhibit a typical trend: error initially decreases with decreasing $R$, reflecting improved local approximation as the nodes become denser. However, beyond a certain threshold, further reduction in spacing leads to increased error in the derivatives due to numerical instability and ill-conditioning in the system matrix.

In the case of MHRBF, the function error remains close to machine precision ($10^{-14}$ to $10^{-12}$) across a broad range of radii ($R \leq 0.5$), before increasing gradually as the spacing becomes too coarse. This flat region followed by error growth reflects both stable convergence and high accuracy. For the derivatives, the error initially decreases with smaller $R$, but then rises again due to amplified round-off error. This behavior is expected, as derivative approximations inherently involve terms inversely proportional to node spacing, which can magnify floating-point errors when $R$ becomes very small. The HRBF method also shows a convergence trend, though with higher overall error.

In summary, both methods exhibit the typical error decay followed by a stability threshold as spacing is varied. However, the MHRBF method demonstrates superior accuracy and smoother convergence profile.

\subsection{Computational Cost Comparison}\label{computationalcost}

\begin{figure}[H]
    \centering
        \begin{subfigure}{0.24\textwidth}
            \centering
            \includegraphics[width=\linewidth]{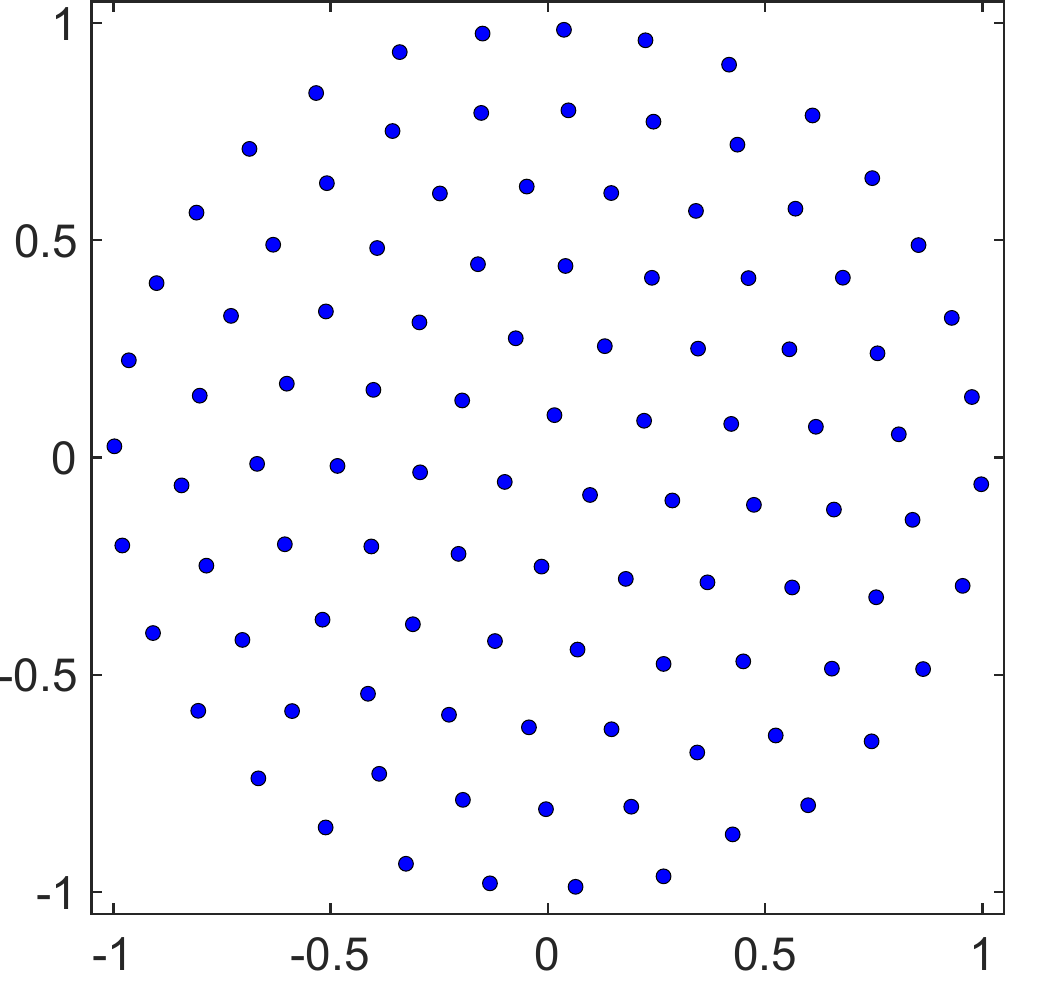}
            \captionsetup{justification=centering}
            \caption{Data nodes}
        \end{subfigure}
        \hspace{0.01\textwidth}
        \begin{subfigure}{0.24\textwidth}
            \centering
            \includegraphics[width=\linewidth]{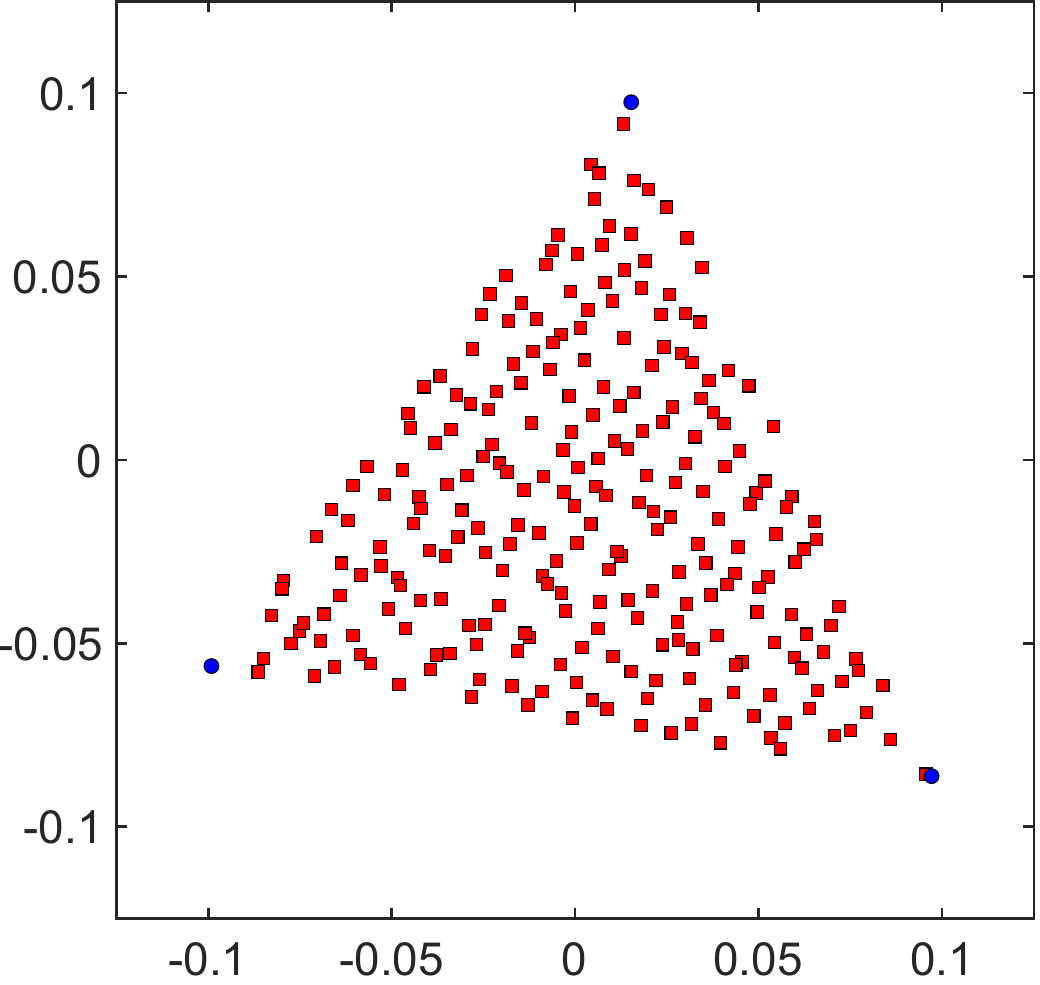}
            \captionsetup{justification=centering}
            \caption{Error evaluation nodes}
        \end{subfigure}
    \captionsetup{font=small}
    \caption{Distribution of data and error evaluation nodes used to evaluate the error as a function of computational cost. The data nodes are blue circles while the error evaluation nodes are red squares.}
    \label{fig:CompCostNodes}
\end{figure}

As a proxy for the computational cost, we evaluate the error against the number of unknowns for a given RBF system. Specifically, the number of unknowns is given by $3N+M$, where $N$ is the number of data nodes and $M=\binom{l+2}{2}$ is the number of polynomial basis functions in 2D for a polynomial of order-$l$. For this analysis the maximum number of data nodes considered is 104, while the error is calculated using 249 evaluation nodes within the triangle of the three data nodes closest to the origin, see \autoref{fig:CompCostNodes}. In all cases, the kernel parameter was set to $n=4$.

\begin{figure}[H]
    \centering
    \begin{subfigure}{0.32\textwidth}
        \includegraphics[width=\linewidth]{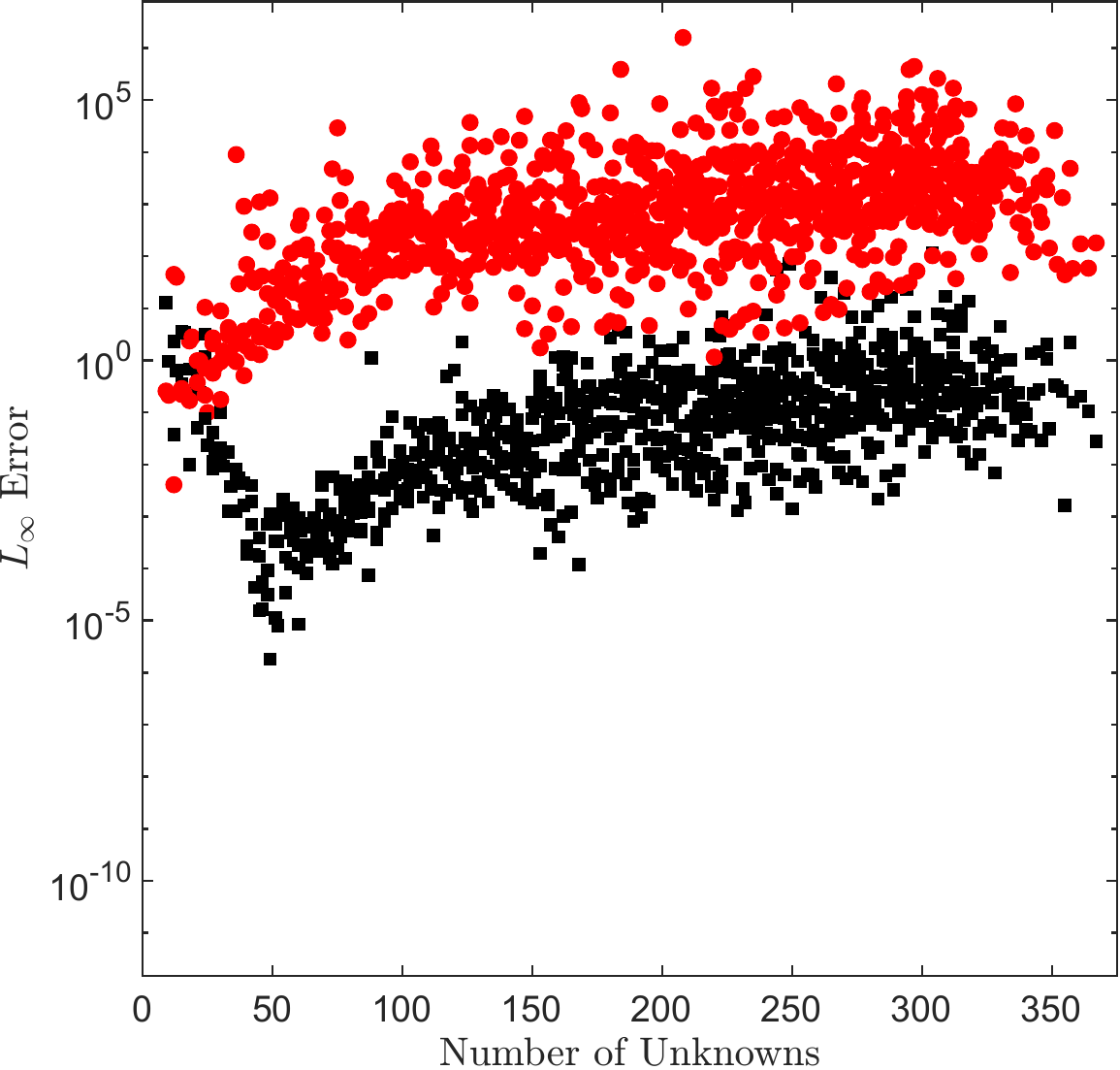}
        \caption{$\varepsilon=0.01$}
    \end{subfigure}
    \hspace{0.05\textwidth}
    \begin{subfigure}{0.32\textwidth}
        \includegraphics[width=\linewidth]{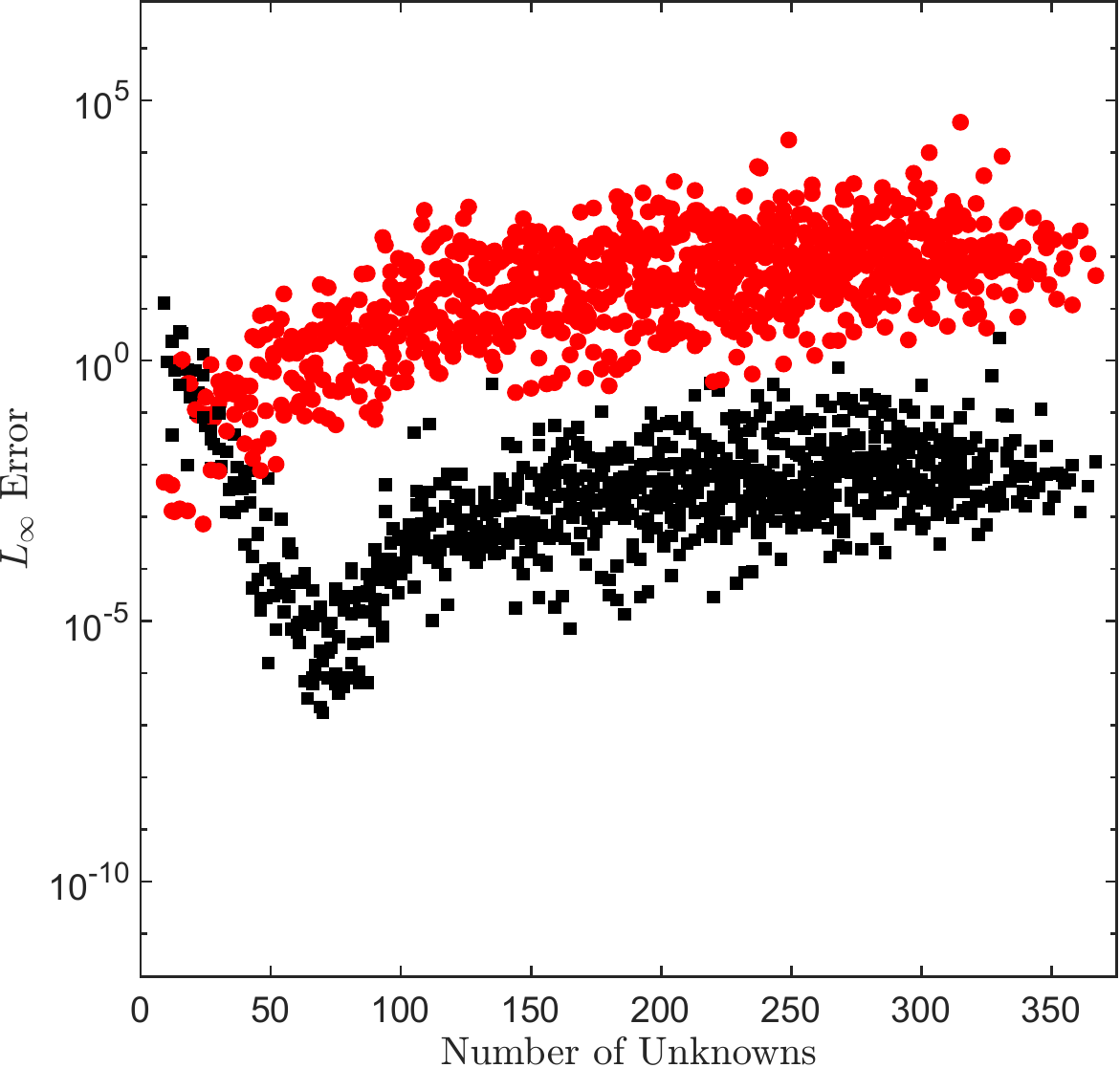}
        \caption{$\varepsilon=0.1$}
    \end{subfigure}
    \\
    \bigskip
    \begin{subfigure}{0.32\textwidth}
        \includegraphics[width=\linewidth]{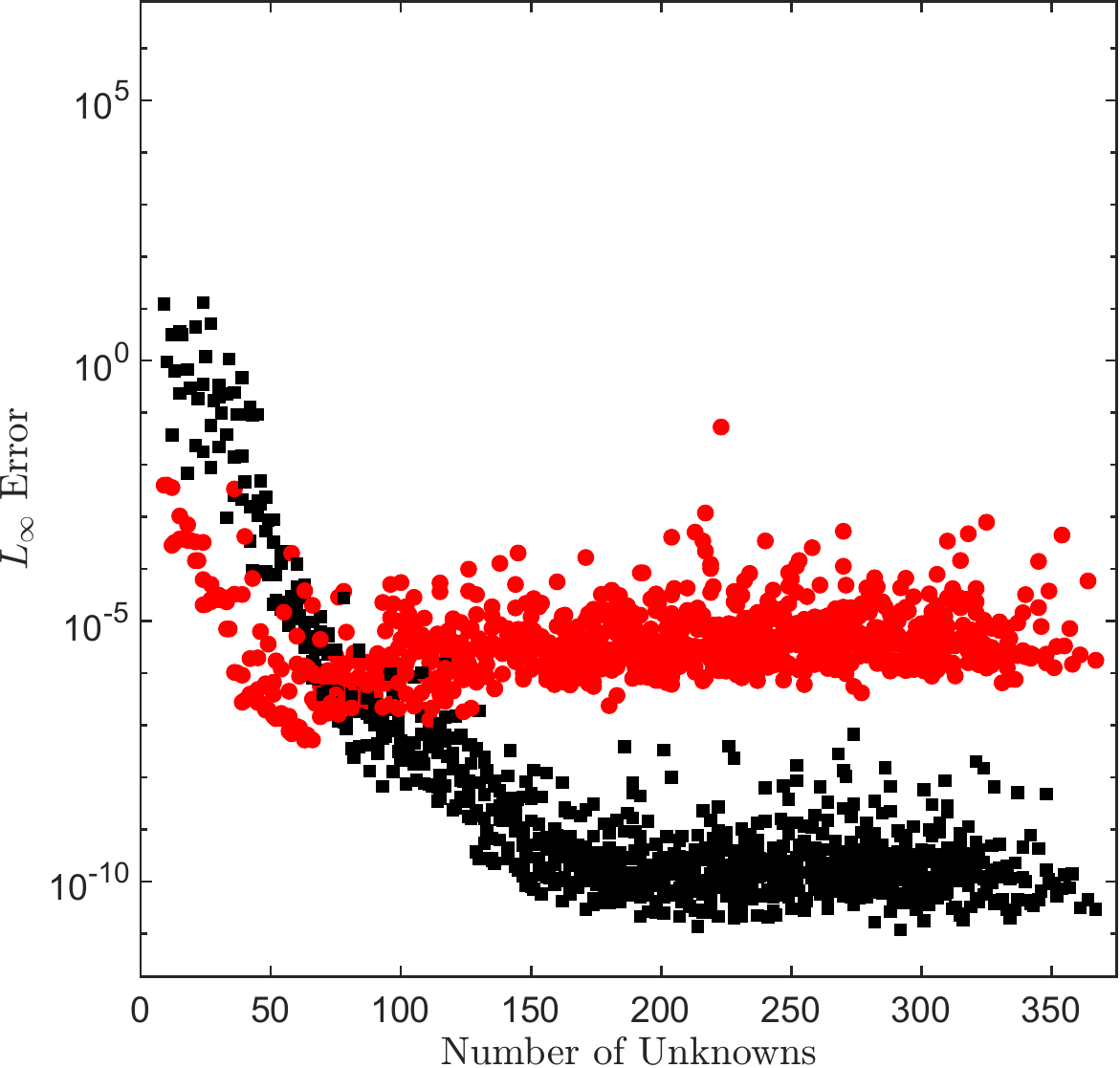}
        \caption{$\varepsilon=1$}
    \end{subfigure}
    \hspace{0.05\textwidth}
    \begin{subfigure}{0.32\textwidth}
        \includegraphics[width=\linewidth]{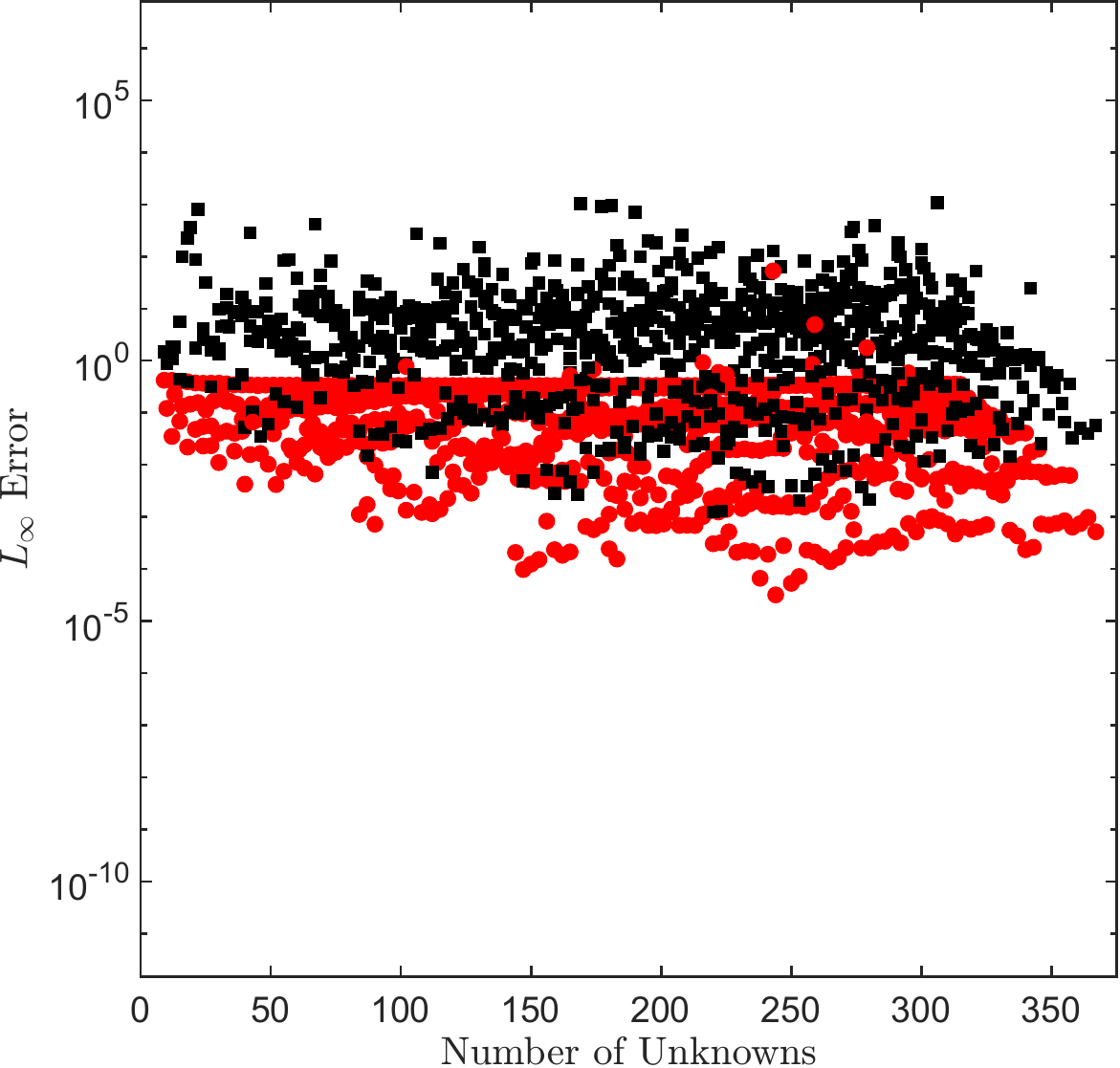}
        \caption{$\varepsilon=10$}
    \end{subfigure}

    \captionsetup{font=footnotesize}
    \caption{Error using the trigonometric test function versus number of unknowns for four shape parameters for data nodes contained in a circle of radius 1 for both the HRBF (red circles) and MHRBF (black squares). No augmentation and augmentation orders 0 to 9 are considered.}
    \label{CompCost_Trig_R1}
\end{figure}

\begin{figure}[H]
    \centering
    \begin{subfigure}{0.32\textwidth}
        \includegraphics[width=\linewidth]{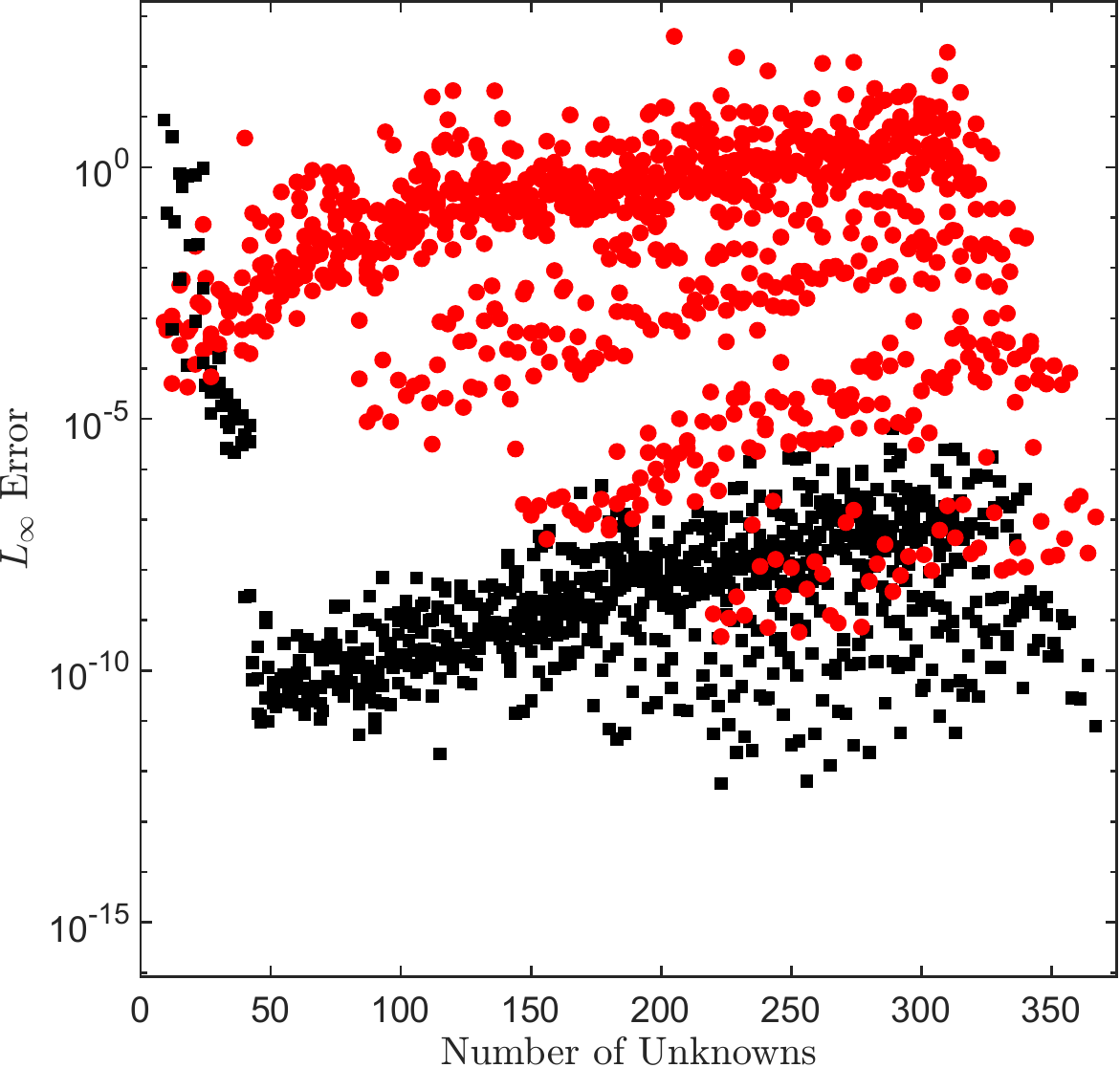}
        \caption{$\varepsilon=0.01$}
    \end{subfigure}
    \hspace{0.05\textwidth}
    \begin{subfigure}{0.32\textwidth}
        \includegraphics[width=\linewidth]{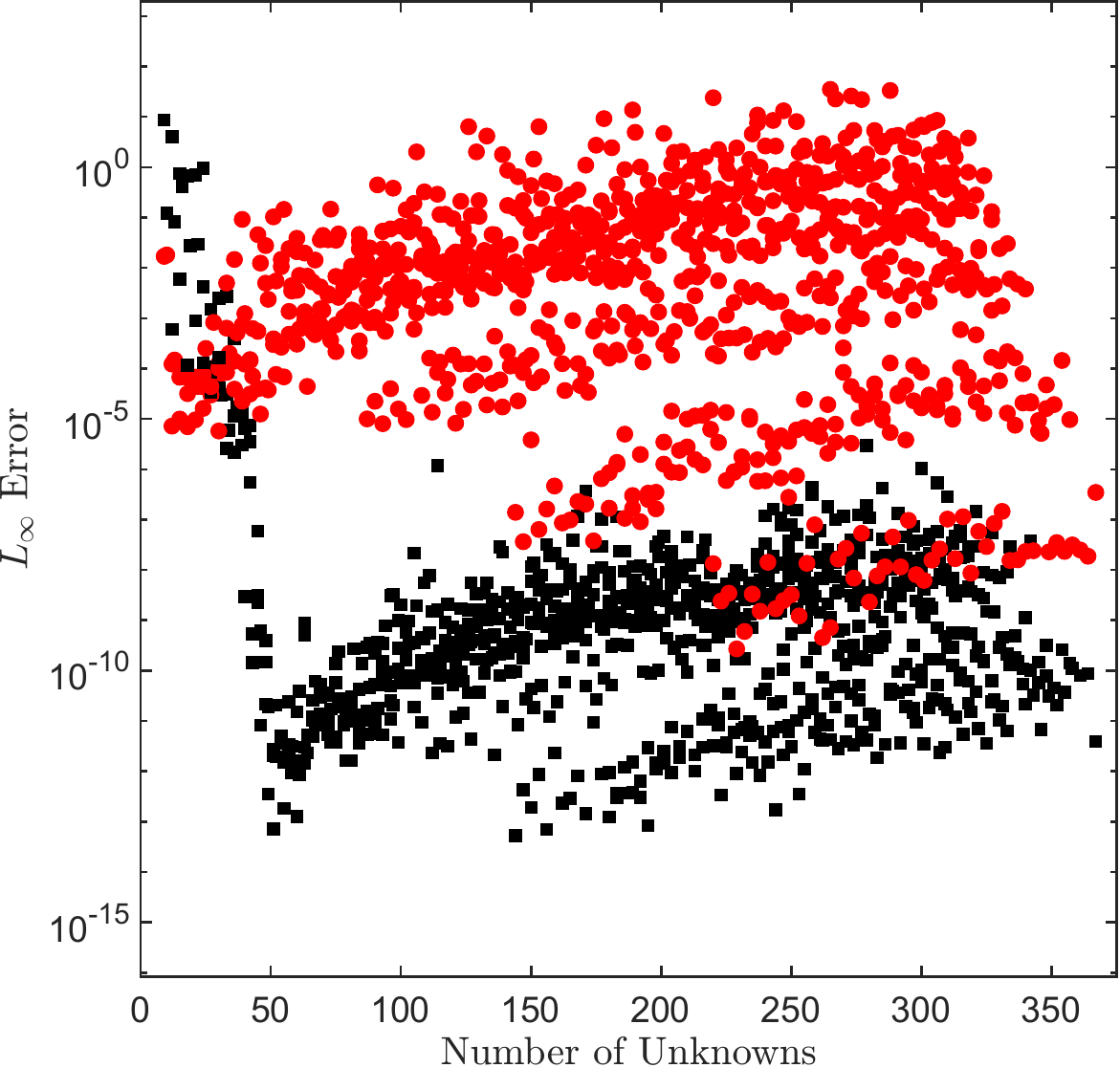}
        \caption{$\varepsilon=0.1$}
    \end{subfigure}
    \\
    \bigskip
    \begin{subfigure}{0.32\textwidth}
        \includegraphics[width=\linewidth]{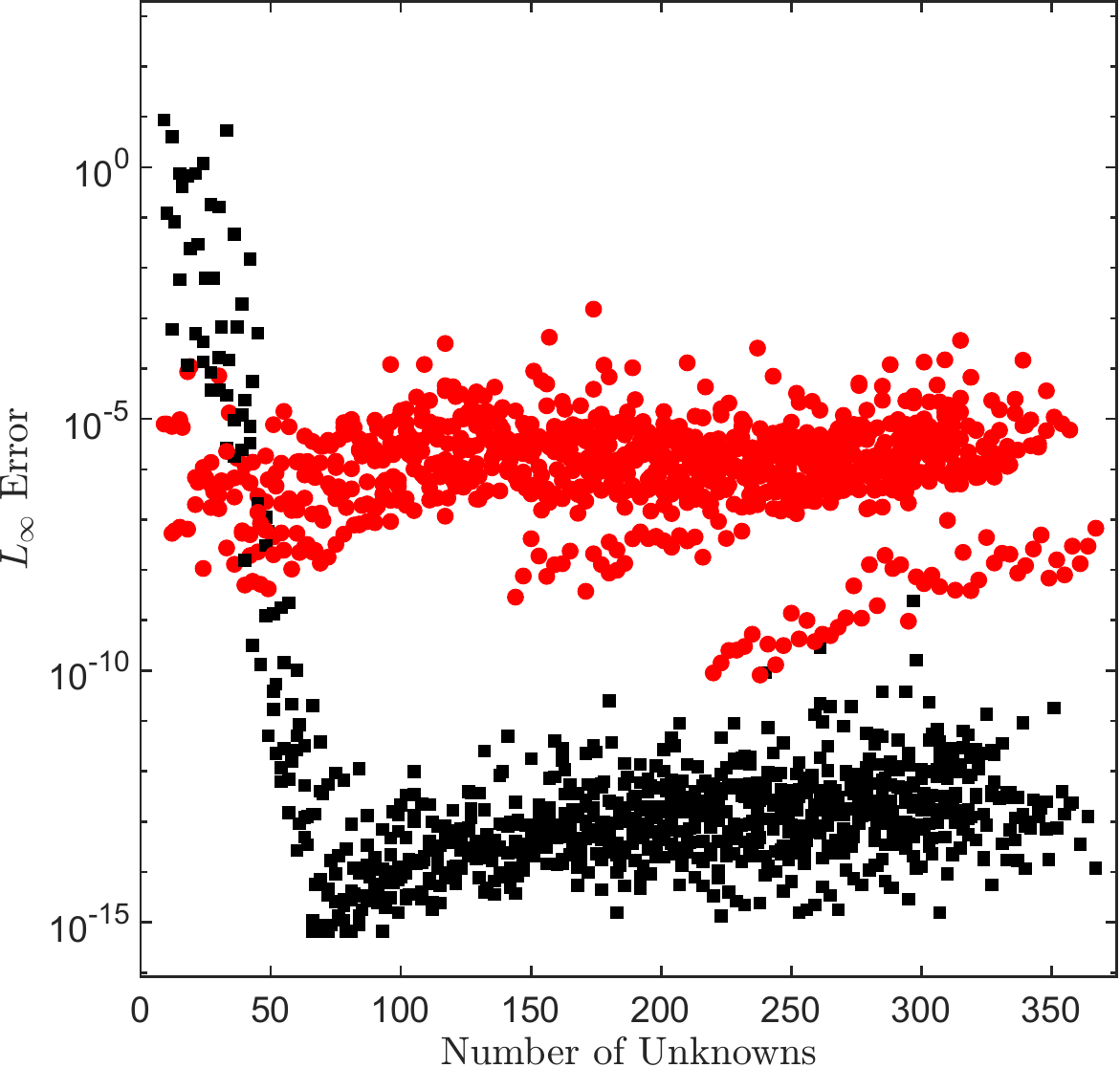}
        \caption{$\varepsilon=1$}
    \end{subfigure}
    \hspace{0.05\textwidth}
    \begin{subfigure}{0.32\textwidth}
        \includegraphics[width=\linewidth]{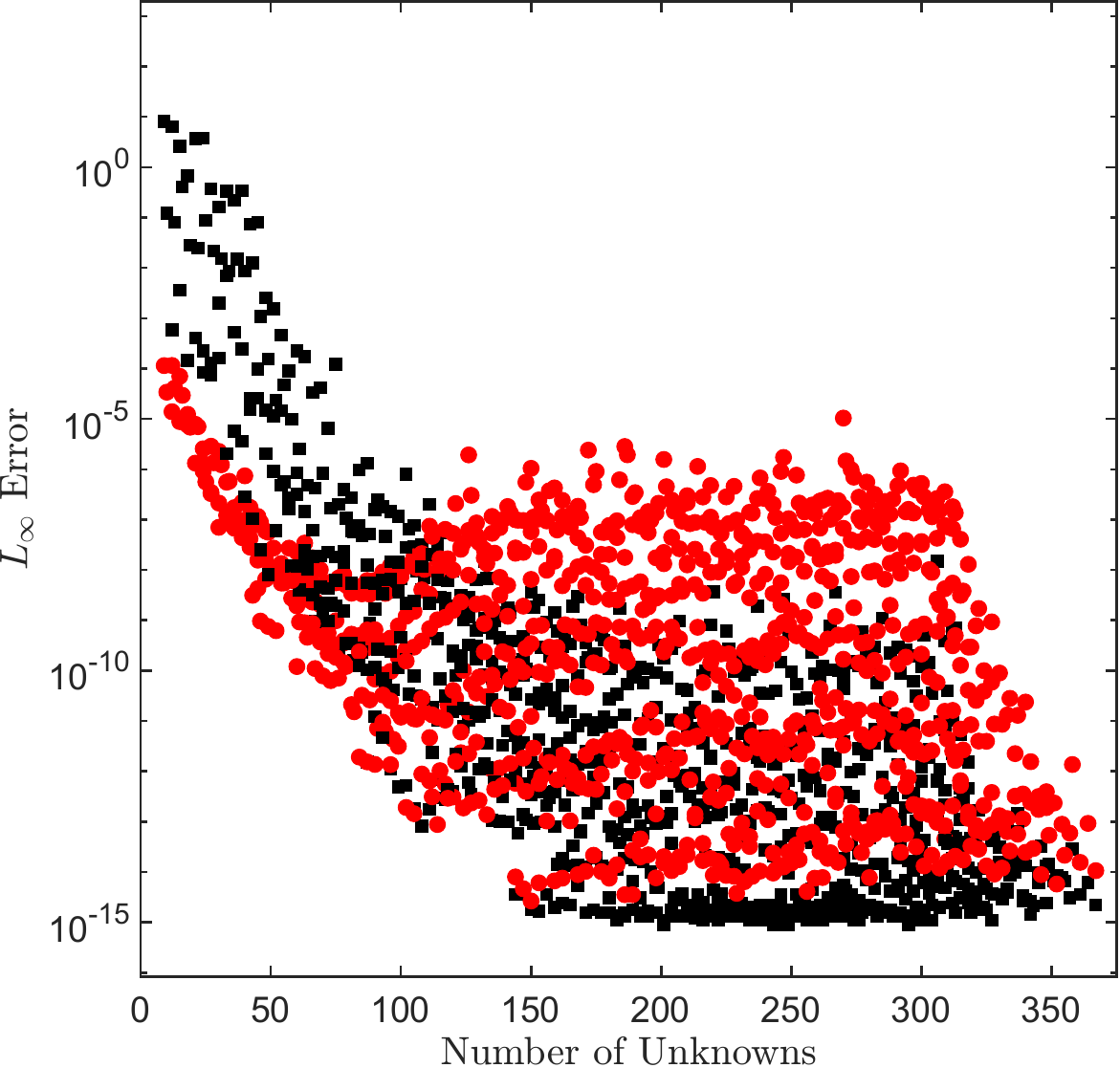}
        \caption{$\varepsilon=10$}
    \end{subfigure}

    \captionsetup{font=footnotesize}
    \caption{Error using the trigonometric test function versus number of unknowns for four shape parameters for data nodes contained in a circle of radius 0.125 for both the HRBF (red circles) and MHRBF (black squares). No augmentation and augmentation orders 0 to 9 are considered. }
    \label{CompCost_Trig_R0p125}
\end{figure}

To begin, we select the three closest data nodes to calculate the HRBF and MHRBF using augmenting polynomial orders ranging from 0 to 9, in addition to no augmenting polynomial. This is done for four shape parameters, $\varepsilon=[0.01, 0.1, 1, 10]$. The error is then computed at the error evaluation nodes, with the largest error recorded. We then repeat this procedure for the four closest nodes, five closest nodes, etc., until all 104 nodes are used. When evaluating the error, we only use polynomial orders which are supported by the number of nodes, i.e. to evaluate an augmenting polynomial order $l$ there must be at least $M$ data nodes. This is performed on data nodes contained in a circle of radius equal to 1 (average spacing of $\approx 2\times 10^{-1}$) and a radius of 0.125 (average spacing of $\approx 2.5\times 10^{-2}$). First, consider the trigonometric test function, \autoref{function1}. The results for a circle radius of 1 is seen in~\autoref{CompCost_Trig_R1} while that for a radius of 0.125 is seen in~\autoref{CompCost_Trig_R0p125}.

There are several important conclusions that can be drawn from the results. With a small number of unknowns, $\lessapprox 40$, the HRBF has a lower error than the MHRBF for the same number of unknowns. For a larger number of unknowns, the MHRBF consistently reaches lower errors at the same cost for $\varepsilon=[0.01, 0.1, 1]$, while the errors are similar using $\varepsilon=10$ and a radius of 0.125. While the errors for HRBF using $\varepsilon=10$ and a radius of 1 are lower for the HRBF, the error values themselves are much higher than what can be achieved using a lower $\varepsilon$. The results also re-iterate the conclusion drawn in \autoref{optimaln2}, which only considered no polynomial augmentation. Specifically, regardless of the average grid spacing or augmenting polynomial order, the choice of $\varepsilon\approx 1$ is the most appropriate value.

\begin{figure}[H]
    \centering
    \begin{subfigure}{0.32\textwidth}
        \includegraphics[width=\linewidth]{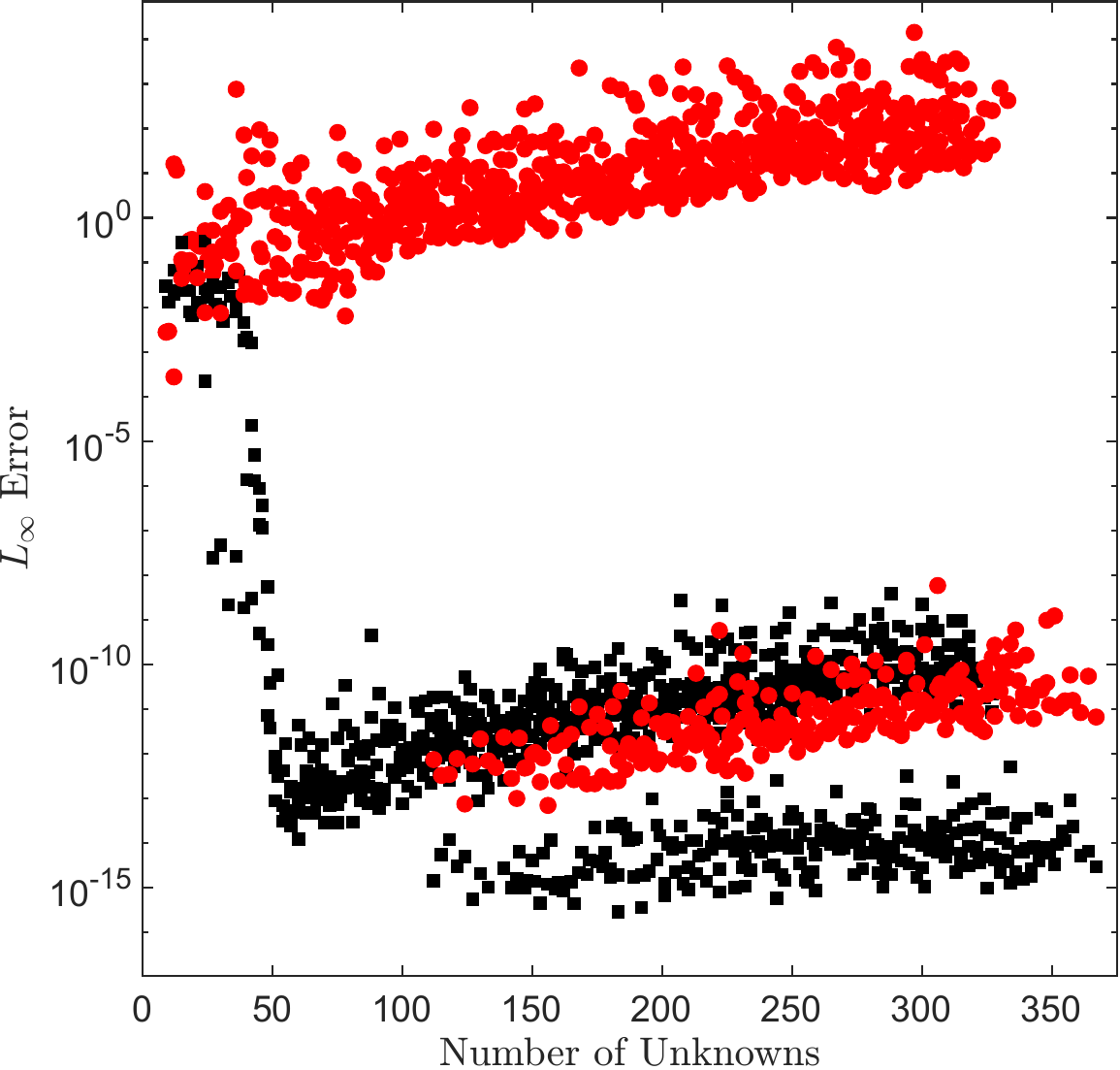}
        \caption{$\varepsilon=0.01$}
    \end{subfigure}
    \hspace{0.05\textwidth}
    \begin{subfigure}{0.32\textwidth}
        \includegraphics[width=\linewidth]{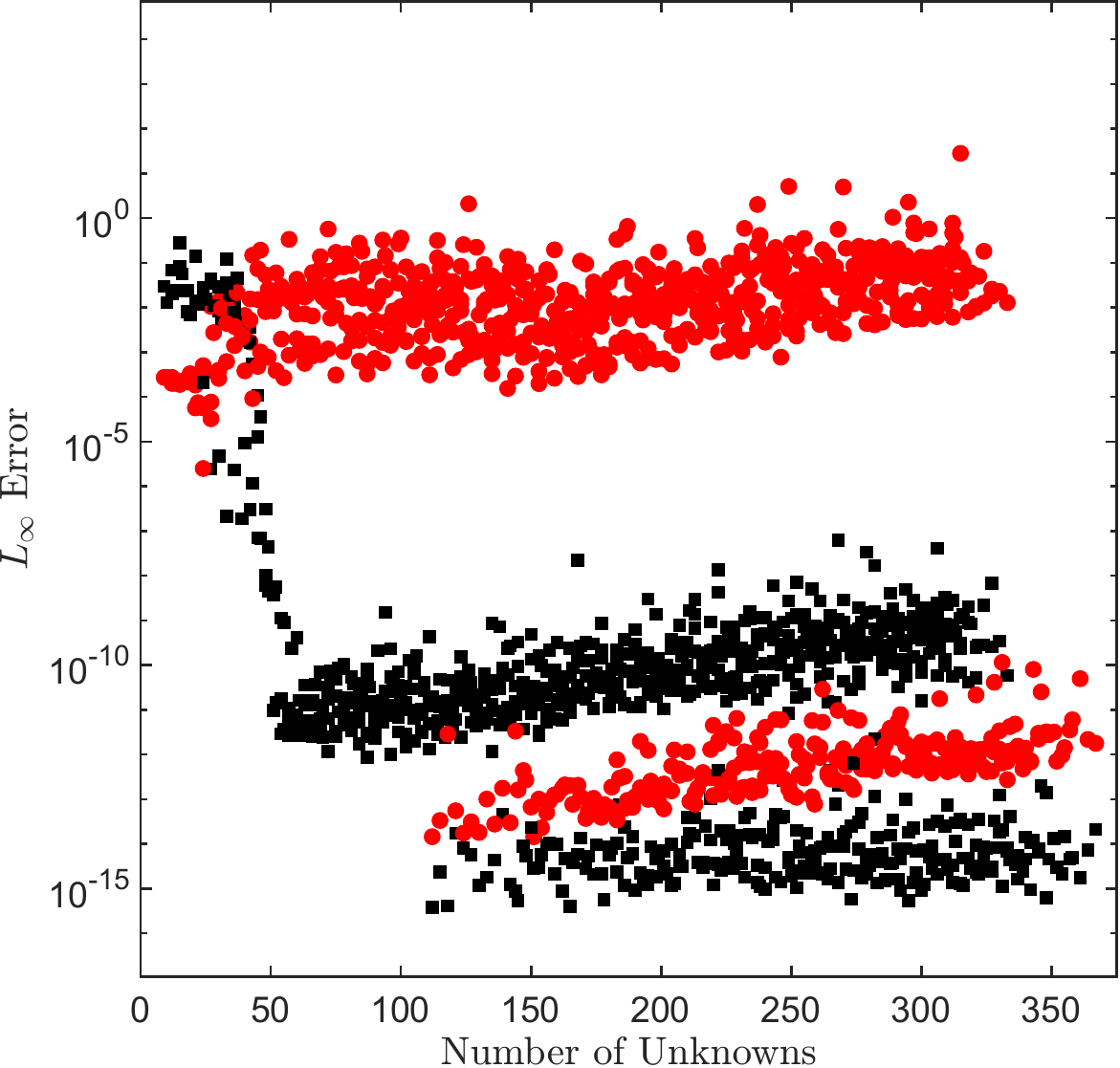}
        \caption{$\varepsilon=0.1$}
    \end{subfigure}
    \\
    \bigskip
    \begin{subfigure}{0.32\textwidth}
        \includegraphics[width=\linewidth]{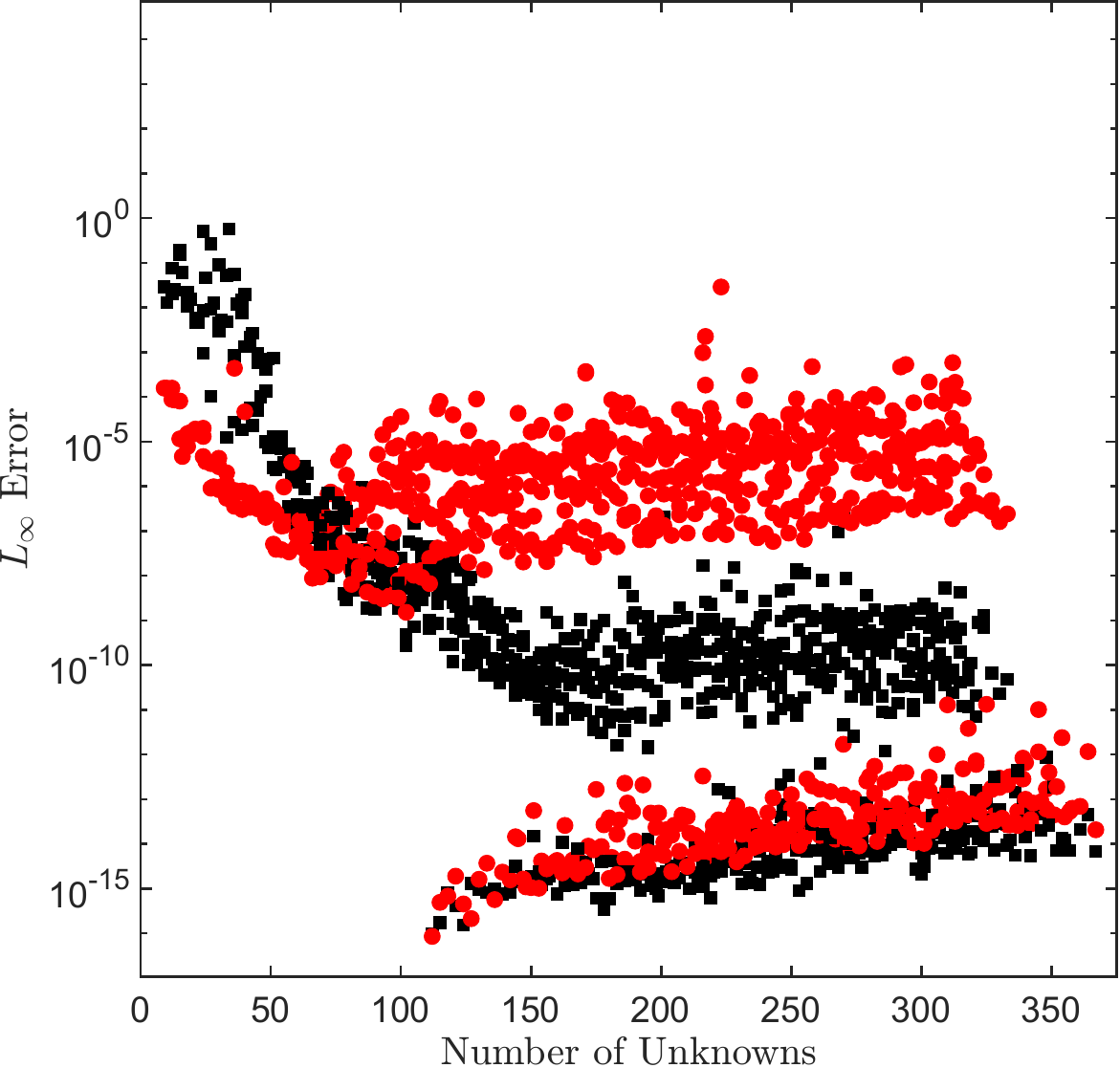}
        \caption{$\varepsilon=1$}
    \end{subfigure}
    \hspace{0.05\textwidth}
    \begin{subfigure}{0.32\textwidth}
        \includegraphics[width=\linewidth]{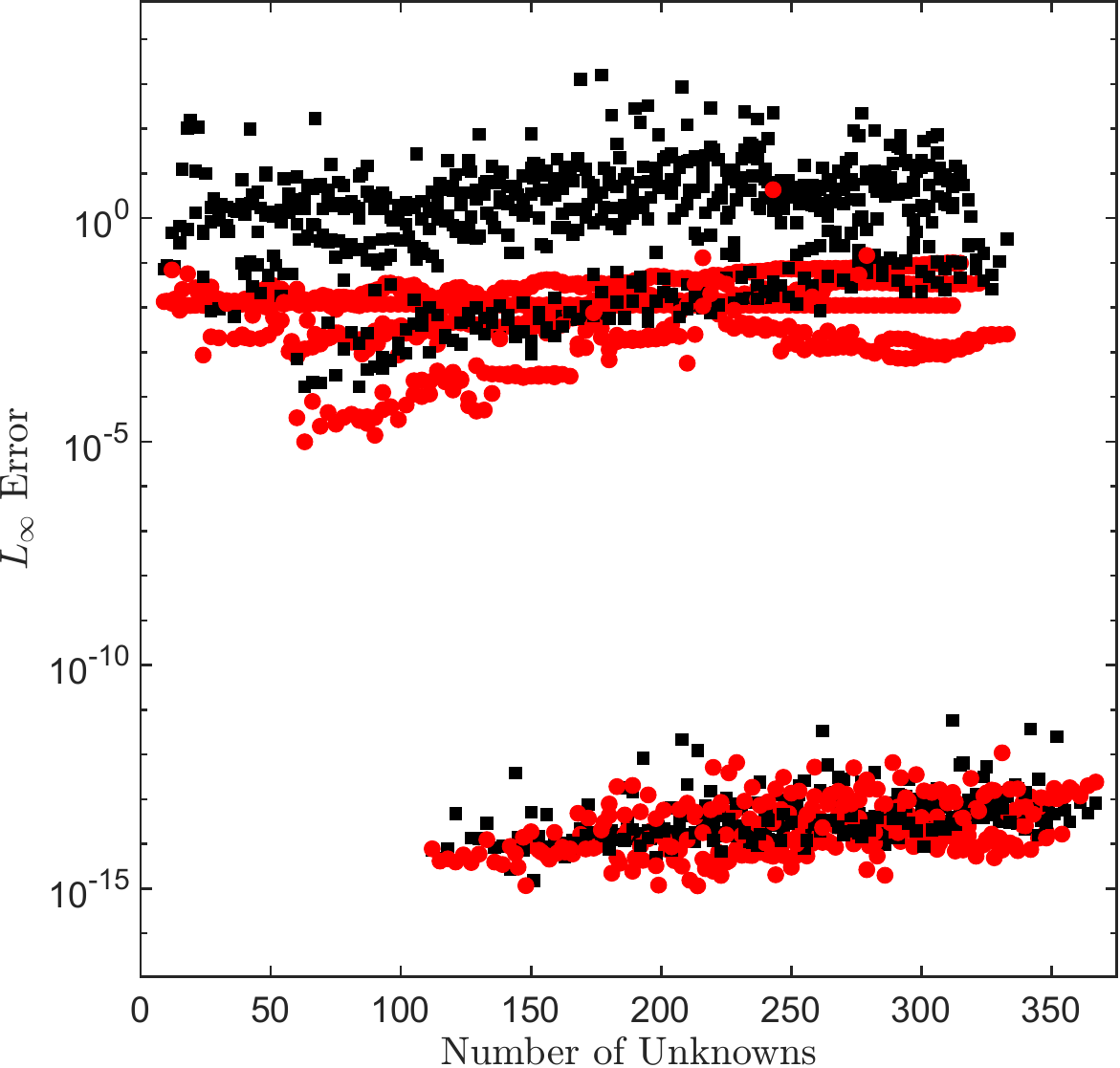}
        \caption{$\varepsilon=10$}
    \end{subfigure}

    \captionsetup{font=footnotesize}
    \caption{Error using the polynomial test function versus number of unknowns for four shape parameters for data nodes contained in a circle of radius 1 for both the HRBF (red circles) and MHRBF (black squares). No augmentation and augmentation orders 0 to 9 are considered. }
    \label{CompCost_Poly_R1}
\end{figure}

\begin{figure}[H]
    \centering
    \begin{subfigure}{0.32\textwidth}
        \includegraphics[width=\linewidth]{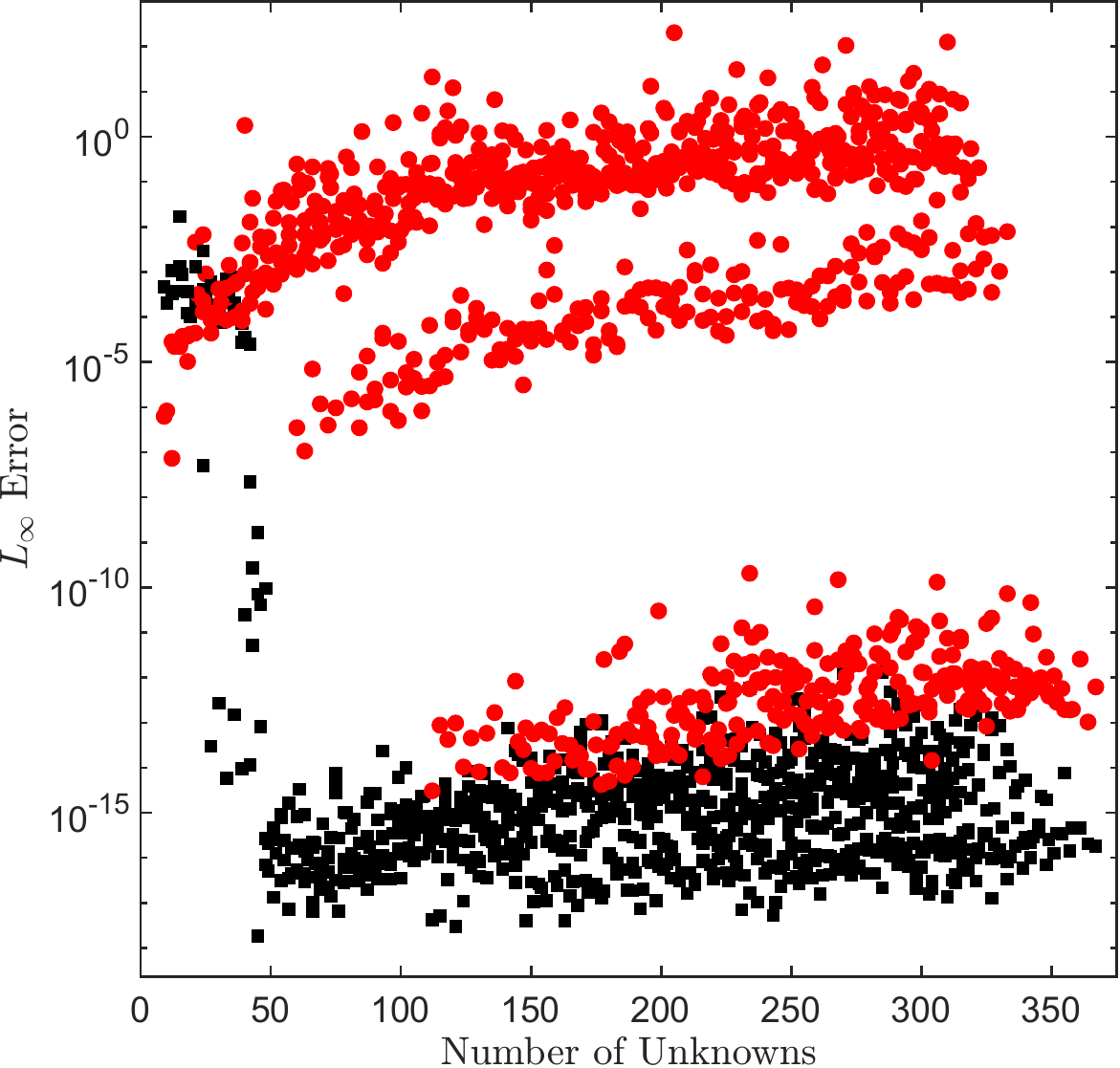}
        \caption{$\varepsilon=0.01$}
    \end{subfigure}
    \hspace{0.05\textwidth}
    \begin{subfigure}{0.32\textwidth}
        \includegraphics[width=\linewidth]{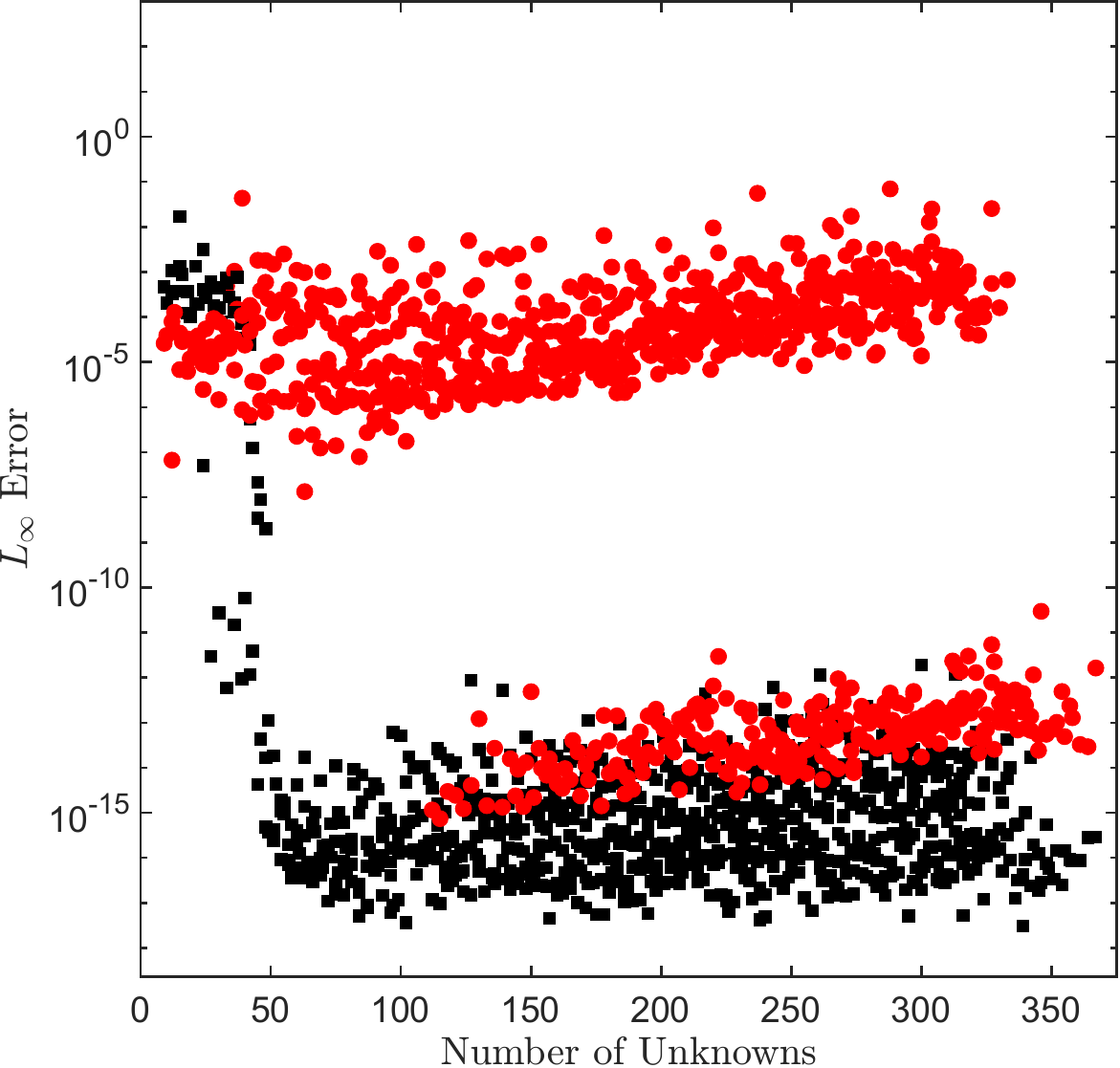}
        \caption{$\varepsilon=0.1$}
    \end{subfigure}
    \\
    \bigskip
    \begin{subfigure}{0.32\textwidth}
        \includegraphics[width=\linewidth]{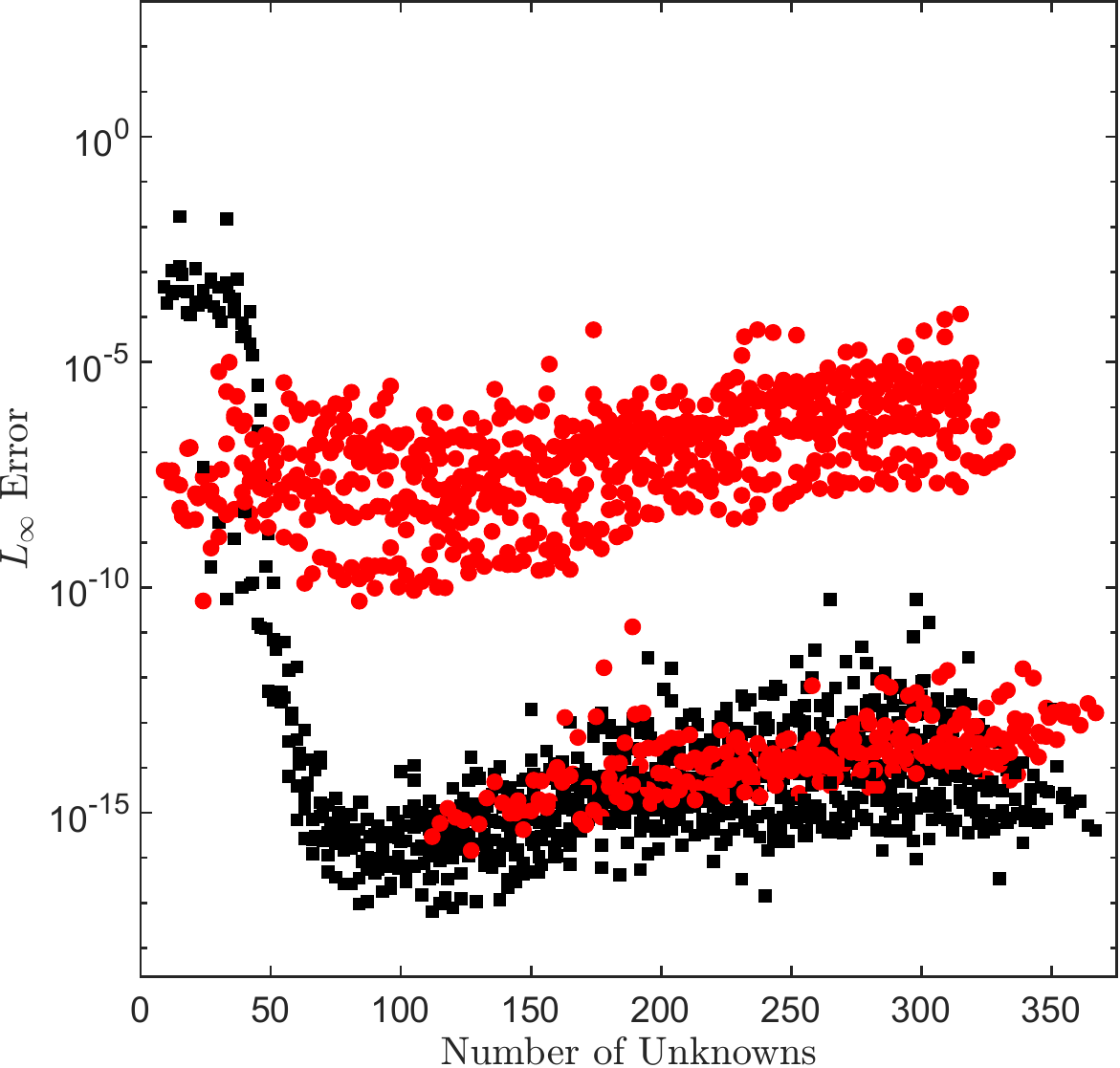}
        \caption{$\varepsilon=1$}
    \end{subfigure}
    \hspace{0.05\textwidth}
    \begin{subfigure}{0.32\textwidth}
        \includegraphics[width=\linewidth]{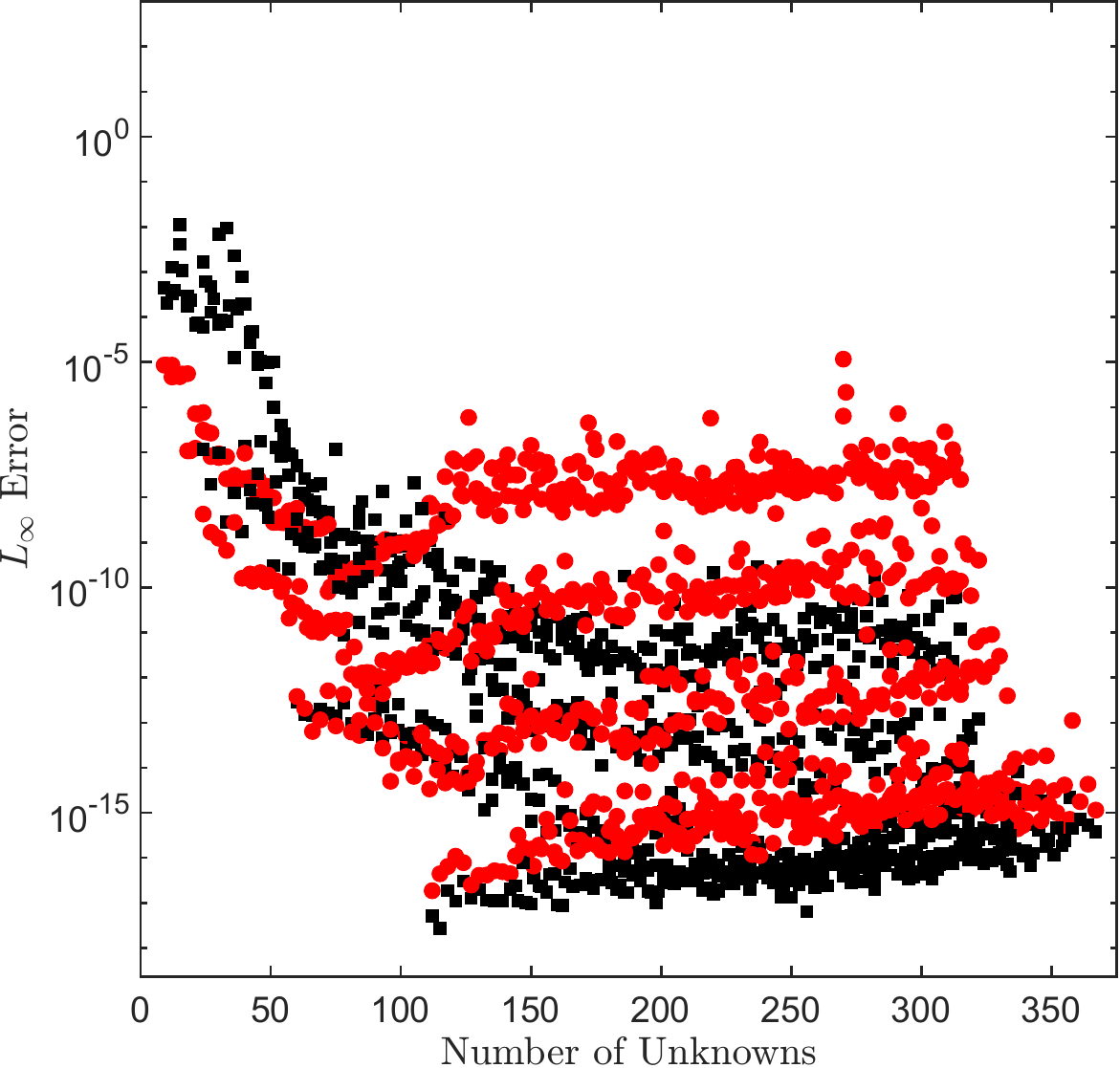}
        \caption{$\varepsilon=10$}
    \end{subfigure}

    \captionsetup{font=footnotesize}
    \caption{Error using the polynomial test function versus number of unknowns for four shape parameters for data nodes contained in a circle of radius 0.125 for both the HRBF (red circles) and MHRBF (black squares). No augmentation and augmentation orders 0 to 9 are considered. }
    \label{CompCost_Poly_R0p125}
\end{figure}

The results of the polynomial test function, \autoref{six-hump}, are shown in \autoref{CompCost_Poly_R1} and \autoref{CompCost_Poly_R0p125}. Many of the same conclusions can be drawn from these results as were observed in the prior function. With a low number of unknowns, the HRBF can outperform the MHRBF in certain situations. Once a certain number of unknowns is met, though, the MHRBF consistently outperforms that HRBF for the same system size. It should be noted that the augmenting polynomial order plays a much larger role in the error for this test function than the previous one, with the error remaining (relatively) large for low augmenting orders and quickly dropping at higher levels. It should be noted that this dependence on the augmenting polynomial is alleviated using the MHRBF, $\varepsilon\lessapprox 1$, and a radius of 0.125, while it is still observed for the HRBF.

\section{Conclusion}

In this work, we introduce the Modified Hermite Radial Basis Function (MHRBF), which uses a modified infinitely-smooth kernel to the Hermite Radial Basis Function (HRBF). Compared to the HRBF, the MHRBF allows for the use of a wider range of the shape parameter while still maintaining accuracy using double-precision mathematics. One important conclusion from this work is that using a kernel monomial parameter of $n\approx 4$ with a shape parameter of $\varepsilon\approx 1$ results in errors close to the minimum possible for two classes of functions, trigonometric and polynomial, regardless of the augmenting polynomial order. This result would allow MHRBF users to select low-order augmenting polynomials, thus reducing the computational cost, and use $n\approx 4$ with $\varepsilon\approx 1$ and remain confident that the interpolation will be close to the smallest possible error given the number of data nodes. Finally, the results also indicate that, with the exception of a low number of nodes, the MHRBF outperforms that HRBF for the same number of unknowns. As the size of the RBF system is directly related to the computational cost and all of the results are obtained using double-precision mathematics, the use of the MHRBF over the HRBF makes this Hermite interpolation technique more attractive for high performance computing situations.

\section*{Acknowledgements}

D.S. is funded by the United States Department of Energy's (DoE) National Nuclear Security Administration (NNSA) under the Predictive Science Academic Alliance Program III (PSAAP III) at the University at Buffalo, under contract number DE-NA0003961.

\bibliographystyle{unsrt}  

\end{document}